%% file: DirichletII2002.tex
\documentclass[a4paper]{amsart}
\usepackage{amssymb}
\usepackage{epsf}

\title
[
Number of facets of three-dimensional Dirichlet stereohedra II
]
{
On the number of facets of three-dimensional Dirichlet stereohedra II:
Non-cubic Groups
}

\author{Daciana Bochi\c{s}}
\author{Francisco Santos}

\address{
Departamento de Matem\'aticas, Estad\'{\i}stica y Computaci\'on,
Universidad de Cantabria, E-39005 Santander, Spain}
\email{dacib@matesco.unican.es, santos@matesco.unican.es}

\date{December 2000, revised April 2002}

\newcommand{\R}{{\mathbb R}}
\newcommand{\Z}{{\mathbb Z}}
\newcommand{\noproof}{\qed}

\newcommand{\Ext}{\hbox{Ext}}
\newcommand{\Infl}{\hbox{Infl}}
\newcommand{\Vor}{\hbox{Vor}}

\newcommand{\Dem}{\noindent {\it Proof: }}

\newtheorem{Teo}{Theorem}[section]

\newtheorem{Prop}[Teo]{Proposition}
\newtheorem{Obs}[Teo]{Remark}
\newtheorem{Exm}[Teo]{Example}
\newtheorem{Cor}[Teo]{Corollary}
\newtheorem{Lema}[Teo]{Lemma}

\begin{document}

\begin{abstract}
We prove that Dirichlet stereohedra for non-cubic crystallographic groups in
dimension 3 cannot have more than 80 facets. The bound depends on the
particular crystallographic group considered and is above 50 only 
on 9 of the 97 affine conjugacy classes of them.

We also construct Dirichlet stereohedra with 32 and 29 facets for a
hexagonal and a tetragonal group, respectively.
\end{abstract}

\maketitle

\section*{Introduction}
This is the second in a series of three papers (see
\cite{MiArt1,MiArt3}) devoted to bound the number of facets that a
Dirichlet stereohedron in 3-space can have.

A \emph{stereohedron} in dimension $d$ is any convex polytope which
tiles the Euclidean space $\R^d$ face-to-face under the action of a
crystallographic group. A \emph{Dirichlet stereohedron} for a
crystallographic group $G$ is the Voronoi
region $\Vor_{GP}(P)$ of a point $P$ in the Voronoi diagram 
of the orbit $GP$.

Finding a good upper bound for the number of facets of $d$-dimensional
stereohedra is mentioned as an important problem in \cite{GruShe-art}
and \cite{Sch-Sen'97}, and it is related to Hilbert's 18$^{th}$
problem \cite{Milnor}. The most relevant previous results on this problem
are:

\begin{itemize}

\item Delone \cite{Delaunay} proved that
a stereohedron of dimension $d$ for a certain
crystallographic group $G$ cannot have more than $2^d (a+1) -2$ facets,
where $a$ is the \emph{number of aspects} of $G$, i.e., the
order of the quotient group of $G$ by its translational subgroup. In
dimension 3, crystallographic groups can have up to 48 aspects,
which gives an upper bound of 390 for the
number of facets of 3-dimensional stereohedra.
Non-cubic groups without reflections have up to 16 aspects (see Table
\ref{Table.Groups}), giving a bound of 134.
\item The 3-dimensional stereohedron with the maximum number of facets
known so far has 38 facets, and it was found by Engel
(see \cite{Engel} and \cite[p. 964]{GruShe-art}). It
is a Dirichlet stereohedron for a cubic group with 24 aspects.
\end{itemize}

Delone's result contrasts with the fact that arbitrarily large
families of non-overlapping congruent convex 3-polytopes 
exist in which the intersection of every two polytopes is a common facet
\cite{Erickson}.

The main results of this series of papers are summarized in Table
\ref{Table.Results}. 
We divide the 219 affine conjugacy classes of 3-dimensional
crystallographic groups in three blocks: Those which contain
reflection planes (100 conjugacy classes) were studied in \cite{MiArt1},
non-cubic groups without reflections (97 classes) are the object of this
paper  and cubic groups without reflections (22 classes) will be
studied in \cite{MiArt3}.
Our methods provide different upper bounds
on the number of facets for each group. The two columns in
Table \ref{Table.Results}
show respectively the global upper bound obtained for each
block of groups and the number of affine conjugacy classes of groups
where our bound is greater than the number of facets of Engel's
stereohedron.

The results for groups with reflections are
specially good, since the upper bound of 18 facets was proved to be
tight in \cite{MiArt1}. For cubic groups we have
indicated the results contained in \cite{MiTesis}, which might be
slightly improved in the final writing of \cite{MiArt3}.
\begin{table}
\begin{center}
\begin{tabular}{|c|c|c|}
\hline
                      & global bound  & nbr. of ``bad'' groups \cr
\hline
Groups with reflections  & & \cr
   \cite{MiArt1}         &     18        &           0      \cr
\hline
Non-cubic groups & & \cr
     $[$This paper$]$    &     84        &          21      \cr
\hline
Cubic groups  & & \cr
   \cite{MiArt3}         &     162       &          18      \cr
\hline
{\bf Total}              &   {\bf 162}   &       {\bf 39}   \cr
\hline
\end{tabular}
\end{center}
\caption{\label{Table.Results} Summary of our results}
\end{table}

Table \ref{Table.numero-maximo-no-cubicos} gives a more detailed
description of the results in this paper. It lists the 21 non-cubic
groups where our bound is greater than 38, together with the specific
bound for each group. The last column in the table indicates where in
this paper the bound is proved. For some of the groups we
have taken Delone's upper bound of $8a+6$.

\begin{table}[htb]
 \renewcommand{\arraystretch}{1.2}
\begin{center}
\begin{tabular}{|c|c|c|cc|}
\hline
Group & Aspects & Planar group & Our upper bound & \dots proved by \\

\hline
 $I\overline{4}c2$    & 8  & $pgg$ & {\bf 40} & Cor. \ref{Coro.pgg.vert} \\
 $P\frac{4_2}{n}\frac{2}{g}\frac{2}{c}$
                      & 16 & $pgg$ & {\bf 40} & Cor. \ref{Coro.pgg.vert} \\
 $R\overline{3}$      & 6  & $p3$  & {\bf 42} & Prop. \ref{Prop.R3}\\
 $R32$                & 6  & $p3$  & {\bf 42} & Prop. \ref{Prop.R3}\\
 $R3c$                & 6  & $p3$  & {\bf 42} & Prop. \ref{Prop.R3}\\
 $I4_1cd$             & 8  & $pgg$ & {\bf 44} & Prop. \ref{Prop.pgg} \\
\hline
 $P3_12$              & 6  & $p1$  & {\bf 48} & Cor. \ref{Coro.planes}\\
 $P3_112$             & 6  & $p1$  & {\bf 48} & Cor. \ref{Coro.planes}\\
 $P6_1$               & 6  & $p1$  & {\bf 48} & Cor. \ref{Coro.planes}\\
 $P4_122$             & 8  & $p1$  & {\bf 50} & Cor. \ref{Coro.P4122}   \\
 $C\frac{2}{a}\frac{2}{c}\frac{2}{c}$
                      & 8  & $p2$  & {\bf 50} & Cor. \ref{Coro.planes}\\
 $I\frac{2}{a}\frac{2}{c}\frac{2}{c}$
                      & 8  & $pgg$ & {\bf 50} & Cor. \ref{Coro.planes}\\
 $P4_12_12$           & 8  & $p1$  & {\bf 64} & Cor. \ref{Coro.planes}  \\
 $I\frac{4_1}{g}$     & 8  & $p2$  & {\bf 70} & Delone\\
 $I4_122$             & 8  & $p2$  & {\bf 70} & Delone  \\
 $I\overline{4}2d$    & 8  & $p2$  & {\bf 70} & Delone  \\
 $F\frac{2}{d}\frac{2}{d}\frac{2}{d}$
                      & 8  & $p2$  & {\bf 70} & Delone \\
\hline
 $P 6_2 2 2$          & 12 & $p2$  & {\bf 78} & Cor. \ref{Coro.planes}\\
 $P 6_1 2 2$          & 12 & $p1$  & {\bf 78} & Prop. \ref{Prop.P6122}\\
 $R \overline{3} \frac{2}{c}$
                      & 12 & $p3$  & {\bf 79} & Cor. \ref{Coro.R32c}\\
 $I\frac{4_1}{g}\frac{2}{c}\frac{2}{d}$
                      & 16 & $pgg$ & {\bf 80} & Cor. \ref{Coro.pgg.vert}\\
\hline
\end{tabular}

\end{center}
 \renewcommand{\arraystretch}{1}
\caption{Non-cubic groups where our upper bound is more 
than 38}
\label{Table.numero-maximo-no-cubicos}
\end{table}

\bigskip
\noindent{\bf Theorem 1}
{\it 
  Dirichlet stereohedra for non-cubic 
  crystallographic groups in dimension 3
   cannot have more than 80 faces. 

  They can
  not have more than 70 except perhaps in the groups
  $R\overline{3}\frac{2}{c}$, $P6_122$, $P6_222$
  and $I\frac{4_1}{g}\frac{2}{c}\frac{2}{d}$.   They can
  not have more than 50 except perhaps in these four and five other groups.
}
\medskip

The upper bounds are proved in the following Sections, starting from more
general methods to more specific ones. For example, Corollary
\ref{Coro.planes} gives an upper bound
derived from simple parameters of each group. This bound is globally
106, and it is higher than 38 in only 34 groups.
The rest of the paper 
concentrates in the top and bottom parts of Table
\ref{Table.numero-maximo-no-cubicos}: lowering the upper bound for
the bad groups in the bottom of the table and removing groups from the
top of the table by showing that they cannot produce more than 38
facets.

Section 4 is essentially
devoted to only one group, the group $P6_122$. We first lower the
upper bound of 96 from Corollary \ref{Coro.planes} to a 78, and then
show the existence of a stereohedron with 32 facets for this group
(Example \ref{Exm.32}). 
This is the highest number of facets 
of a stereohedron for a non-cubic group
obtained to date. 

This example shows, on the one hand, that the general
upper bound of 80 is at most 2.5 times worse than the actual maximum. And on
the other hand, that non-cubic groups can produce stereohedra almost as
complicated as cubic groups.
It is worth noting that our construction method shows why it is
natural to expect stereohedra with ``many'' facets for the group
$P6_122$. In fact, Lemma \ref{Lema.screw} implies that any
group whose name (in the International Crystallographic Notation, see
\cite{Lockwood}) contains the string $i_12$, produces Dirichlet
stereohedra with at least $4i+1$ facets.
We show a stereohedron with 29 facets for the
group $I4_122$ (Example \ref{Exm.29}).

\section{A first upper bound for each crystallographic group}
\label{Sec.Firstbound}
All throughout the paper $G$ will denote a 3-dimensional
crystallographic group without reflections and $P\in \R^3$ will be a
base point for an orbit $GP$, so that the Dirichlet stereohedron we
want to study is the (closed) Voronoi region $\Vor_{GP}(P)$. 
We assume
that $P$ has trivial stabilizer under the action of $G$, since
otherwise $GP$ is also an orbit of a proper subgroup of $G$.
Further, there is no loss of generality in assuming
both $P$ and the metric parameters of $G$ (the parameters which
identify $G$ in its affine conjugacy class of groups) to be
sufficiently generic since a small
perturbation to either $P$ or $G$ can only increase the number of
facets in the Dirichlet stereohedron induced.
We say that a point $Q\in GP$ is a neighbor of
$P$ if $\Vor_{GP}(P)$ and $\Vor_{GP}(Q)$ share a facet. Hence, the
number of facets of $\Vor_{GP}(P)$ equals the number of neighbors of
$P$.

\subsection{Outline of our method}
\label{Sub.Outline}

Let $v\in\R^3$ be a vector such that the corresponding
translation $\tau_v$ is in $G$. We will call vertical the direction of
$v$ and horizontal the planes orthogonal to it.
Let $Z_P$ denote the closed infinite band of width
$2|v|$ bounded by the horizontal planes passing through
$\tau_vP$ and $\tau_{-v}P$.

\begin{Lema}
\label{Lema.ZP}
All the neighbors of $P$
lie in $Z_P$. Moreover, the only points in the boundary
of $Z_P$ which can be neighbors of $P$ are $\tau_vP$ and
$\tau_{-v}P$.
\end{Lema}

\Dem Let $Q\in GP$. Suppose that $Q$ is outside $Z_P$ or that it is on
its boundary but it equals neither $\tau_vP$ nor $\tau_{-v}P$.
Without loss of generality assume that $PQ\cdot v$ is positive
(otherwise change $v$ to $-v$).  Consider the points $P'=\tau_vP$ and
$Q'=\tau_{-v}Q$. The four points $PP'QQ'$, in this order, are the
vertices of a planar quadrilateral whose angles at $P'$ and $Q'$ are
greater or equal than 90 degrees. Since the four points are in $GP$,
$P$ and $Q$ cannot be neighbors.  \qed

In $Z_P$ there are only a finite family $\alpha_1,\dots,\alpha_k$ of
horizontal planes containing points of $GP$. Our method consists on
bounding the number of neighbors of $P$ in each $\alpha_i$ separately
and using the sum of the numbers obtained as a bound for the total
number of neighbors. In all cases of interest to us the 
direction will be chosen so that 
the subgroup $G_0$ of $G$ consisting of 
horizontal motions (motions which send every
horizontal plane to itself) contains two independent
translations. That is to say, it is a 2-dimensional crystallographic group
acting on each horizontal plane.

Since we are assuming $P$ to be sufficiently generic, all the
intersections $GP\cap \alpha_i$, $i\in\{1,\dots,k\}$ are orbits of
$G_0$. The following statement is Lemma 1.3 in \cite{MiArt1}.

\begin{Lema}
\label{Lema.twoplanes}
Under these assumptions, given two different
horizontal planes $\alpha_i$ and $\alpha_j$, a necessary condition for
$Q\in \alpha_i\cap GP= G_0Q$ to be a neighbor of $P\in \alpha_j\cap
GP=G_0P$ is that the two 2-dimensional Dirichlet stereohedra
$\Vor_{G_0P}(P)$ and $\Vor_{G_0Q}(Q)$ overlap (we are implicitly
projecting $G_0P$ and $G_0Q$ along the direction of $v$ to a common
horizontal plane).
\noproof
\end{Lema}

If, moreover, every
element of $G$ sends horizontal planes to horizontal planes,
then we have the following result.
We recall that the normalizer of a group $G_0$ in a bigger
group $H$ is the subgroup of $H$ consisting of elements $g\in H$ such
that $g^{-1} G_0 g=G_0$.

\begin{Lema}
\label{Lema.normalizer}
Let $G$ be a 3-dimensional crystallographic group. Let a vertical
direction be chosen so that every element of $G$ sends horizontal
planes to horizontal planes. Let $G_0$ be the subgroup of
horizontal motions of $G$.

Then the orthogonal projection of any orbit of $G$ to a horizontal
plane is contained in an orbit of the normalizer of $G_0$ in the group
of affine isometries of that plane.
\end{Lema}

\Dem Let $g\in G$ and let $h\in G_0$. Since $g$ sends horizontal
planes to horizontal planes and $h$ sends horizontal planes to
themselves, $g^{-1}hg$ sends horizontal planes to themselves. Hence
$g^{-1}hg\in G_0$ (i.e. $G_0$ is normal in $G$).  This, together with
the fact that $g$ composed with the projection to any horizontal plane
is an affine isometry on that plane, implies the statement.  \qed

The following result, most of which comes from \cite{MiArt1},
gives a first upper bound on
the number of neighbors of $P$ in each $\alpha_i$, depending on the
type of the group $G_0$:

\begin{Teo}[\protect{\cite[Theorem 3.1]{MiArt1}}]
\label{Teo.dos.orbitas} Let $G_0$ be a planar crystallographic group
and let $G_0P$ and $G_0Q$ be two orbits of points with trivial
stabilizer. Then, the number of Dirichlet regions of one of the orbits
overlapped by each Dirichlet region of the other orbit is bounded
above by:
\begin{itemize}

\item[(i)] one, for groups of types $pmm$, $p3m1$, $p4m$ and $p6m$.
\item[(ii)] two, for groups of types $cm$, $cmm$, $p31m$ and $p4g$.
\item[(iii)] four, for groups of types $p1$, $p3$, $p4$, $p6$, $pmg$ and $pm$.
\item[(iv)] seven, for the groups $p2$ and $pg$.
\item[(v)] eleven, for the group $pgg$.
\item[(vi)] seven, for the group $pgg$, if $G_0P\cup G_0Q$ is contained in
an orbit of the normalizer of $G_0$ in the affine isometries of $\R^2$
(see Proposition \ref{Prop.normal.pgg}).
\end{itemize}
The bounds are all tight except perhaps those for $pg$ and $pgg$.
\end{Teo}
\Dem
Parts (i) to (iv) appear in \cite{MiArt1} as Theorem 3.1.
Part (v) was also proved in  \cite{MiArt1}, although not explicitly stated
as a result (see Remark 3.2 and the paragraph before Lemma 3.5 in that
paper). Part
(vi) is proved as Proposition \ref{Prop.normal.pgg} in Section 3
of this paper. \qed

We conjecture that the upper
bound for $pgg$ is seven even without the normalizer condition.

\subsection{The classification of 3-dimensional
crystallographic groups}
\label{Sub.Classification}

The first two invariants to affinely classify a crystallographic group
$G$ are its \emph{translational subgroup} and its \emph{point group}.
The translational subgroup $T$ is the subgroup consisting of
translations, and belongs to one of the 14 \emph{Bravais lattice types}.
The point group of $G$ is the quotient group $[G:T]$. The
elements of $[G:T]$ are the \emph{aspects} of $G$.
The point group can be understood as
a discrete group of motions in the 2-dimensional unit
sphere satisfying the ``crystallographic restriction'': there are no
rotations of order 5 or greater than 6. There are 32 groups satisfying
this, modulo isometries of the sphere, and all of them appear as point
groups of some crystallographic group.

In turn, the lattices are classified into seven
\emph{crystallographic systems}:
\emph{monoclinic},
\emph{triclinic}, \emph{orthorhombic}, \emph{hexagonal},
\emph{trigonal}, \emph{tetragonal} and  \emph{cubic}. This
classification is according to the point groups of their normalizers in
the group of affine isometries of $\R^3$, except
that most authors include in the trigonal system some
groups whose lattice is hexagonal, the
distinction between the hexagonal and trigonal systems relying on
properties of the point group.

Groups which contain reflection planes were
studied in \cite{MiArt1} and cubic groups will be dealt with in
\cite{MiArt3}. Table \ref{Table.Groups} lists the 97 affine conjugacy classes
of non-cubic 3-dimensional crystallographic groups without
reflections, divided according to crystallographic systems and point
groups. 
We use the International Crystallographic Notation to name 
2-dimensional and 
3-dimensional crystallographic groups and have taken the monograph 
\cite{Lockwood} as a source for their classification. 
The capital letter at the
beginning of a name indicates the Bravais type of the lattice within
its crystallographic system. The rest of the
name encodes generators for the group in a certain way.

\input{table97.tex}

Only 58 of the 97 groups in Table \ref{Table.Groups} have more than 4 aspects.
We will deal only with these 58 groups since the Delone
bound for groups with four or less aspects is already 38, matched
by Engel's example. Moreover, experimental evidence indicates that in
groups with few aspects
Delone's bound is a good approximation to what can happen
in the worst case: There are Dirichlet stereohedra for groups with one
and two aspects having 14 and 20 facets respectively
\cite[p. 963]{GruShe-art}.

The 58 groups belong to the following crystallographic systems:
orthorhombic (10 groups), trigonal (13 groups), hexagonal (9 groups)
and tetragonal (26 groups).
They are listed in Tables \ref{Table.Orthorhombic}--\ref{Table.Tetragonal}.

\input{table58.tex}

There are the following 8 possible Bravais types of lattices
in the 58 groups of Tables
\ref{Table.Orthorhombic}--\ref{Table.Tetragonal}.
\begin{itemize}
\item The three primitive lattices (noted $P$ in Table
\ref{Table.Groups}), generated by a vertical translation and the
following planar lattices in horizontal planes: a rectangular lattice
in the orthorhombic system, a triangular lattice in the hexagonal and
trigonal systems, and a square lattice in the tetragonal system (see
Figure \ref{Fig.lattices}).
Each primitive lattice has a primitive cell associated to it, which is
an orthogonal prism over a rectangle, triangle or square
respectively. The primitive cell tiles the space face-to-face
producing as vertices of the tiling an orbit of the lattice.

\begin{figure}[htbp]
\centerline{
  \epsfysize=2 cm
  \epsfbox{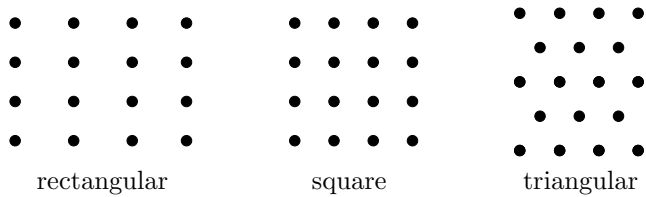}
}
\centerline{\null \hskip .5 cm
rectangular \hskip 1.8 cm square \hskip 1.7 cm triangular
}
  \caption{The three planar lattices}
  \label{Fig.lattices}
\end{figure}

\item The two body-centered lattices (noted $I$) in the orthorhombic
and tetragonal systems. They 
are generated by the primitive lattice and any translation sending
a vertex of the primitive cell to the centroid of the cell. Hence,
they contain the primitive lattice as a sublattice of index 2.

\item The face-centered lattice in the orthorhombic system (noted $F$),
generated by the primitive lattice and translations from a vertex of
the primitive cell to the centroids of the three facets of the cell
incident to that vertex. It 
contains the primitive lattice as a sublattice of index 4.

\item The base-centered lattice in the orthorhombic system (noted $C$),
generated by the primitive lattice and the translation from a vertex of
the primitive cell to the centroid of one of the three facets of the cell
incident to that vertex. It 
contains the primitive lattice as a sublattice of index 2.

\item The rhombohedral lattice in the trigonal system (noted $R$),
generated by the
primitive lattice and a translation from one of the vertices of the
primitive cell (a triangular prism) to a point in the axis of the
primitive cell and whose distances to the two triangular bases of the
prism are in the ratio 1:2. It 
contains the primitive lattice as a sublattice of index 4.
\end{itemize}

\subsection{A first upper bound}
\label{Sub.Firstbound}

Suppose that $G$ is one of the 58 groups of Tables
\ref{Table.Orthorhombic}--\ref{Table.Tetragonal}. The lattice of $G$ 
contains as a sublattice one of the three primitive
lattices, and we choose as vertical the direction of the principal
axis in the corresponding
primitive cell. The principal axis is well defined in the
hexagonal, trigonal and tetragonal systems. In the orthorhombic system
there are three perpendicular and equivalent
axes. We choose one arbitrarily
unless the lattice is base-centered (noted $C$) in which case we
choose as vertical the direction perpendicular to the
centered faces.

This choice fulfills the conditions
required in Lemmas \ref{Lema.twoplanes} and \ref{Lema.normalizer}:
$G_0$ is crystallographic and
every element of $G$ sends horizontal planes to horizontal planes.  
We take as vector $v$ to apply Lemma \ref{Lema.ZP} the shortest
vertical translation in $G$.  The following Lemma shows how to compute
the number of horizontal planes in the band $Z_P$ of Lemma
\ref{Lema.ZP} which contain points of $GP$.
This number depends only on the type of the lattice and the 
number of aspects of $G$ and $G_0$.

\begin{Lema}
\label{Lema.planes}
Apart of the two boundary planes and the plane containing $P$, there
are exactly $2al/a_0-2$ horizontal planes in $Z_P$ containing points
of  the orbit $GP$, where:
\begin{itemize}
\item $a$ is the number of aspects of $G$.
\item $a_0$ is the number of aspects of the horizontal subgroup $G_0$.
\item $l$ is a number depending only on the lattice of $G$, and takes
the value 1 in primitive and base-centered lattices, 2 in
body-centered and face-centered lattices, and 3 in rhombohedral lattices.
\end{itemize}
\end{Lema}
\Dem Let $T$ and $T_0$ denote the translational subgroups of $G$ and
$G_0$, respectively. The number $l$ is the
ratio of horizontal planes containing points of $T$ to
planes containing points of the primitive lattice, so that each
orbit of $T$ projects to $l$ orbits of $T_0$.

The closed band $Z_P$ contains $2al$ plus $a_0$ orbits of $T_0$.
Each horizontal plane containing points of $GP$ contains
$a_0$ orbits of $T_0$. Hence, there are $2al/a_0+1$ such planes,
including the two boundary planes and the one containing
$P$. Subtracting 3 gives the statement.  
\qed

\begin{Cor}
\label{Coro.planes}
Let $i$ be the maximum number of intersections between two
different orbits of the planar crystallographic group $G_0$, as stated
in Theorem \ref{Teo.dos.orbitas}. Then, the number of
neighbors of the point $P$ in an orbit $GP$ is bounded above by
\[
2i(al/a_0-1)+8,
\]
where $a$, $a_0$ and $l$ are as in Lemma \ref{Lema.planes}.
\end{Cor}

\Dem
By Theorem \ref{Teo.dos.orbitas}, $P$ has at most $i$ neighbors in
each horizontal plane of $Z_P$ other than the two boundary ones and
the one containing $P$. The only possible neighbors in the boundary
are $\tau_v(P)$ and $\tau_{-v}(P)$, and in the plane containing $P$
there are at most six neighbors, since a planar Dirichlet stereohedron
has at most six facets.
\qed

In Tables \ref{Table.Orthorhombic}--\ref{Table.Tetragonal}, the column
``Cor. \ref{Coro.planes}'' indicates the upper bound given by this
statement in each of the 58 groups of interest to us. 

\section{Influence region. Groups with a horizontal $p2$ or $p3$}
\label{sec.influence}

Corollary \ref{Coro.planes} is based on using global bounds for the
number of regions of $\Vor_{G_0Q}$ which can intersect the region
$\Vor_{G_0P}(P)$. But of course the number of regions intersected
depends on $Q$, and we are interested only in a few concrete
possibilities for the point $Q$ in each 3-dimensional crystallographic
group. To study the dependence of the number of intersections with $Q$
we recall the formalism of {\em extended Voronoi regions} and {\em
influence regions}, already used in \cite{MiArt1}.

\subsection{Influence regions in planar crystallographic groups}
Let $G_0$ be a crystallographic group in the plane.
Let $N_0$ be a crystallographic group containing $G_0$.
We will call \emph{fundamental subdomains} of
$G_0$ (w.r.t. $N_0$) the fundamental domains of $N_0$.
For a fixed fundamental subdomain $D$ of $G_0$, we call
\emph{extended Dirichlet region} of $D$ under $G_0$ (and denote it
$\Ext_{G_0}(D)$) any region containing $\cup_{Q\in D} \Vor_{G_0Q}(Q)$.

Observe that, for any given $D$, the fundamental subdomains $\{\tau D:
\tau\in N_0\}$ tile the plane.
We call \emph{influence region} of $D$ under $G_0$ (and denote it
$\Infl_{G_0}(D)$) the union of all the tiles $\tau D$ with $\tau\in N_0$
for which $\Ext_{G_0}(\tau D)$
and $\Ext_{G_0}(D)$ overlap.

\begin{Lema}
Let $P$ and $Q$ be two points in the plane and let $D$ be a fundamental
subdomain containing $P$. A necessary condition for $\Vor_{G_0P}(P)$
and $\Vor_{G_0Q}(Q)$ to overlap is that $Q$ lies in $\Infl_{G_0}(D)$.
\end{Lema}

\Dem
Let $\tau\in N_0$ be such that $Q\in \tau D$. By construction,
$\Vor_{G_0P}(P)\subset \Ext_{G_0}(D)$ and
$\Vor_{G_0Q}(Q)\subset \Ext_{G_0}(\tau D)$. This proves the statement.
\qed

Observe finally that the action of $G_0$ divides the fundamental
subdomains in a finite number of orbits. Then, for each point $Q\in
\R^2$, the number of regions of $\Vor_{G_0Q}$ which intersect
$\Vor_{G_0P}(P)$ can be bounded above by counting how many fundamental
subdomains in the orbit of fundamental subdomains containing $Q$ intersect
the influence region
$\Infl_{G_0}(D)$.

In this section we apply this formalism to 
2-dimensional crystallographic groups of types 
$p2$ and $p3$:

\subsubsection*{Groups of type $p3$:}
Let $G_0$ be a group of type $p3$, generated by an order
3 rotation and two translations of equal length with angle
$\pi/6$ to one another.
The order-3 rotation centers of the group form a triangular
lattice, and any two adjacent elementary triangles of the lattice form
a fundamental domain. Our choice of $N_0$ is generated by
$G_0$ together with any reflection on one of the sides of an
elementary triangle. Each elementary triangle is a fundamental
subdomain, and there are two $G_0$-orbits of fundamental subdomains,
which we call $B$ and $W$ for ``black'' and ``white''. The three
neighbors of a black triangle are white and vice-versa. We assume the
subdomain containing $P$ to be white.

\begin{figure}[htbp]

\centerline{
  \epsfxsize=10cm
  \epsfbox{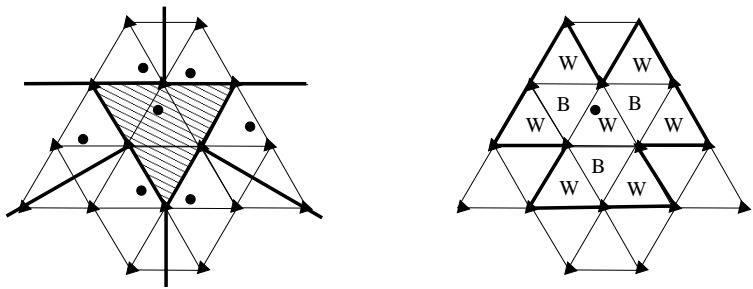}
}
  \caption{The influence region for a group of type $p3$}
  \label{Fig.p3}
\end{figure}

The extended Dirichlet region of an elementary triangle consists of
this triangle and the three adjacent to it. This is shown in the left
part of Figure \ref{Fig.p3} in the following way: a base point $P$ in
a fundamental subdomain $D$ and
other six points of the orbit $G_0P$ are drawn. Let $\tau P$ be one of the
six points. The region of the plane consisting of points closer to
$\tau P$ than to $P$ for any choice of $P$ in the central fundamental
subdomain is clearly disjoint with $\Ext_{G_0}(D)$. What is left is
precisely the union of $D$ and its three adjacent subdomains.

Hence, the influence
region consists of the fundamental subdomains which are either in 
$\Ext_{G_0}(D)$ or adjacent to a
fundamental subdomain in $\Ext_{G_0}(D)$. This gives the 10 fundamental
subdomains in the right part of Figure \ref{Fig.p3}.

In total, there are seven white
and three black triangles in the influence region.

\begin{Lema}
\label{Lema.p3}
Let $G_0$ be a group of type $p3$. Let $P$ and $Q$ be two points in
the plane and suppose that $P$ lies in a white triangle, in the
sense explained above. Then, the number of regions of $\Vor_{G_0Q}$
which overlap $\Vor{G_0P}(P)$ is at most three if $Q$ is in a black
triangle and at most four if it is in a white triangle.
\end{Lema}

\Dem
The bound of three for black triangles is obvious since there are only
three black triangles in the influence region. For the white triangles
we would in principle have a bound of seven, but part (iv) of Theorem
\ref{Teo.dos.orbitas} gives four.
\qed

\subsubsection*{Groups of type $p2$ with rectangular grid:}
Let $G_0$ be a group of type $p2$, generated by two independent
translations and a rotation of order two. We assume further that the
two generating translations are orthogonal to one another, which
implies that the order-2 rotation centers form a rectangular grid. Any
two adjacent elementary rectangles in this grid form a fundamental
domain for $G_0$. Our choice of $N_0$ is generated by
$G_0$ together with any reflection on one of the sides of an
elementary rectangle. Then, each elementary rectangle is a fundamental
subdomain and there are two $G_0$-orbits of fundamental subdomains
which we call $B$ and $W$ for ``black'' and ``white''. The four
neighbors of a black rectangle are white and vice-versa. We assume the
rectangle containing $P$ to be white.

\begin{figure}[htbp]
\centerline{
  \epsfxsize=12.5cm
  \epsfbox{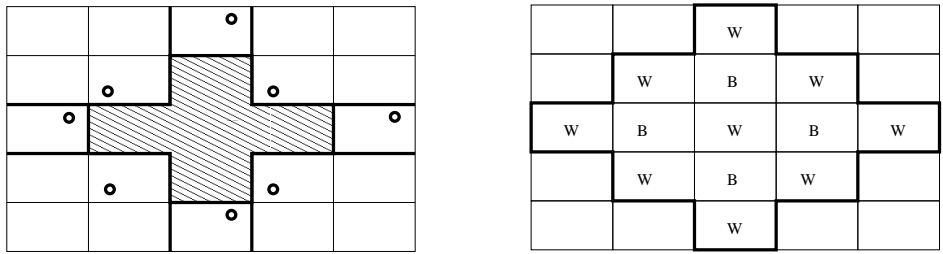}
}
  \caption{The influence region for a group of type $p2$}
  \label{Fig.p2}
\end{figure}

The extended Dirichlet region of an elementary rectangle consists of
this rectangle and the four adjacent to it. The influence
region consists of the extended Dirichlet region and the eight
elementary rectangles adjacent to it. See Figure \ref{Fig.p2}.
In total, there are nine white
and four black rectangles in the influence region.

\begin{Lema}
\label{Lema.rectangular.p2}
Let $G_0$ be a group of type $p2$ with a rectangular grid. Let $P$ and $Q$ be
two points in the plane and suppose that $P$ lies in a white rectangle, with
the meaning explained above. Then, the number of regions of $\Vor_{G_0Q}$
which overlap $\Vor{G_0P}(P)$ is at most four if $Q$ is in a black
rectangle and at most seven if it is in a white rectangle.
\end{Lema}

\Dem
The bound of four for black rectangles comes from the influence
region. The bound of seven for white rectangles comes from Theorem
\ref{Teo.dos.orbitas}.
\qed

We are interested in the case where $G_0$ is the horizontal subgroup of a
3-dimensional crystallographic group $G$. We will always make sure
that $N_0$ is a
group of horizontal motions that not only contains $G_0$ but is also
contained in the normalizer of $G$ in the isometry group of $\R^3$.
The following lemma implies that if this happens then there is no
loss of generality in assuming that our base point $P$ lies
over any fixed fundamental subdomain $D$, since if this is not the case then
there
will be an orbit $GP'$ isometric to $GP$ (in particular, with congruent
Dirichlet stereohedra) and with $P'$ over $D$.

\begin{Lema}
Let $G$ be a 3-dimensional crystallographic group. Let $G_0$ be its
horizontal subgroup and let $N_0$ be a horizontal group containing
$G_0$ and contained in the normalizer of $G$.
Let $D$ be a fundamental subdomain of $G_0$. Let $P$ be any point in
$\R^3$. Then, there is an isometry $\tau$
sending $P$ to a point $P'$ over $R$ and with $\tau(GP)=G(\tau P)$.
\end{Lema}

\Dem Let $\tau\in N_0$ be a horizontal motion sending $P$ to a point
$P'$ over $D$. This exists since $D$ is a fundamental domain of
$N_0$. Since $\tau$ is in the normalizer of $G$, $\tau GP=G\tau
P=GP'$. This finishes the proof.
\qed

\subsection{Groups with a horizontal $p2$ of rectangular type}

\begin{figure}[htb]
 \epsfysize=1.35cm
  \centerline{\epsfbox{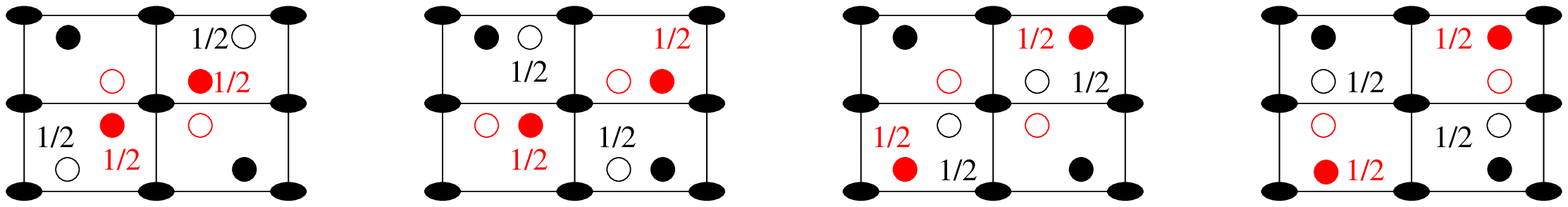}}
{\centerline{
$P\frac{2}{n}\frac{2}{n}\frac{2}{n}$ \hspace{1.4cm}
$P\frac{2}{a}\frac{2}{n}\frac{2_1}{n}$ \hspace{1.4cm}
$P\frac{2}{n}\frac{2_1}{c}\frac{2_1}{c}$ \hspace{1.4cm}
$P\frac{2}{a}\frac{2_1}{c}\frac{2}{c}$
}}
\vspace{0.5cm}
 \epsfysize=1.95cm
  \centerline{\epsfbox{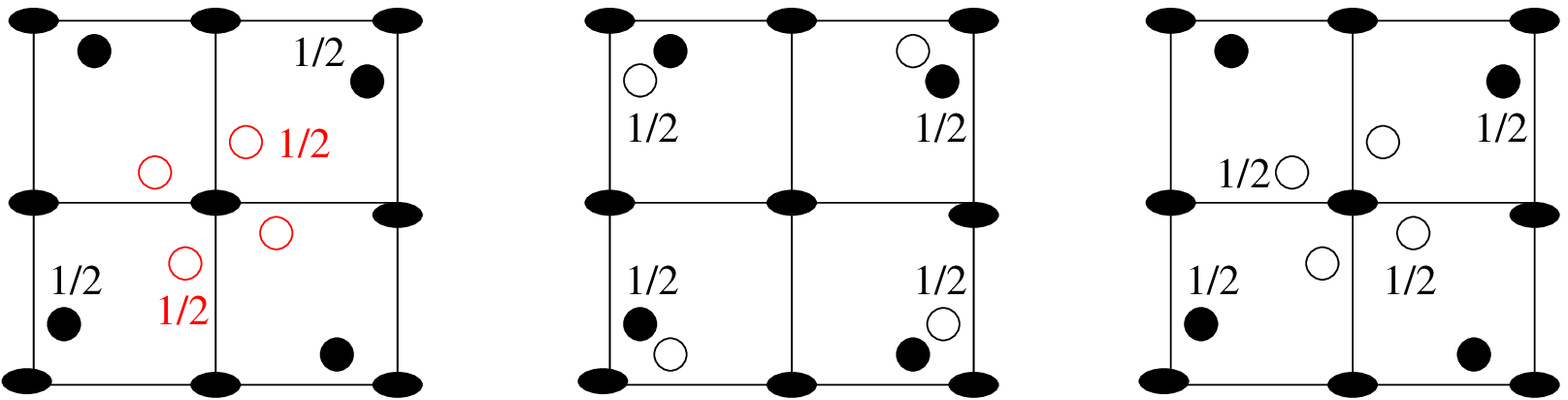}}
  {\centerline{
\null\hspace{.2 cm}
$P\frac{4_2}{n}$  \hspace{1.8cm}
$P4_222$ \hspace{1.6cm}
$P4_22_12$ }}
\vspace{0.5cm}
 \epsfysize=1.95cm
  \centerline{\epsfbox{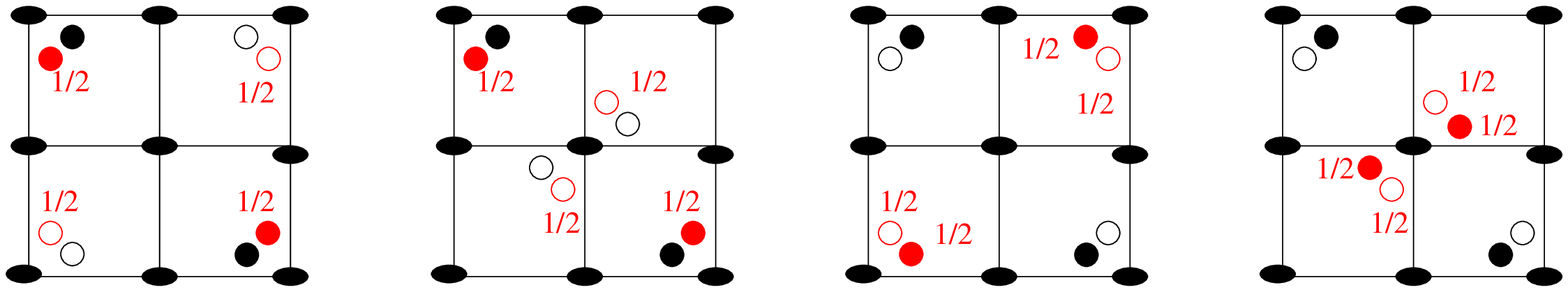}}
  {\centerline{
$P\overline{4}2c$  \hspace{1.8cm}
$P\overline{4}2_1c$  \hspace{1.5cm}
$P\overline{4}c2$  \hspace{1.6cm}
$P\overline{4}n2$ }}
  \caption{Orthorhombic and tetragonal groups whose horizontal group is a
rectangular $p2$}
\label{Fig.rectangular.p2}
\end{figure}

\begin{Prop}
\label{Prop.rectangular.p2} The four orthorhombic and seven
tetragonal groups with 8 aspects displayed in Figure
\ref{Fig.rectangular.p2} have a horizontal group $G_0$ of type $p2$
with rectangular grid. The reflections on vertical planes containing
edges of elementary rectangles of the grid lie in the normalizer of
$G$.

In each case, the band $Z_P$ contains six horizontal planes (not
counting the boundary and middle ones) of a generic orbit $GP$, four
which produce black orbits and two which produce white orbits. Hence,
their Dirichlet stereohedra cannot have more than 38 facets.
\end{Prop}

\Dem The bound of 38 follows from the rest of the statement by a counting
argument, using Lemma \ref{Lema.rectangular.p2}: to the 6 neighbors in the
plane of $P$ and 2 in the boundary of $Z_P$ we have to add four in each of
the four ``black'' planes of orbits and seven in the two ``white'' ones.

The rest of the statement can be easily checked in the graphical
representation of each
group displayed
in Figure \ref{Fig.rectangular.p2}. Let us first explain this representation.

In each picture, a horizontal projection of a generic orbit appears. The
black ovals in the corners of the rectangles/squares represent the order-2
rotation centers of the horizontal group $p2$. Hence, any of the four
rectangles is a fundamental subdomain, any adjacent two are a fundamental
domain of $G_0$, and the four together are a fundamental domain of the
translational subgroup of $G_0$. 
In each picture, the black circle in the
bottom-right quadrant represents the base point $P$ for the 3-d orbit $GP$.
The other seven circles represent points in the other seven aspects of
$G$, with their heights indicated in the following way: let $z$ be the
height of the base point 
and assume that the shortest vertical translation in
$G$ has length 1. Then, black circles with a number $\alpha$ represent
horizontal orbits at heights $z+\alpha$ (and hence at $z+\alpha+m$ for any
integer $m$) while white circles with a number $\alpha$ represent orbits at
height $-z+\alpha$ (and hence at $-z+\alpha +m$ for any integer $m$). The
absence of a number means $\alpha=0$ (as happens in the base point itself).

To make all this clearer, let us describe generators for the group
$P\frac{2}{n}\frac{2}{n}\frac{2}{n}$ deduced from its graphical
representation. First, there is the translational subgroup, generated by
horizontal translations on the sides of the big rectangle and a vertical
translation of length 1. Then, we have to describe how to get 
the other seven points in the figure from the base point. 
The white point in the
bottom-right rectangle is obtained by a rotation of order two around the
vertical axis of the bottom-right rectangle, followed by a reflection on the
horizontal plane at height zero.
The two points on the upper-left rectangle are obtained from the
two in the bottom-right rectangle by a rotation of order two on the vertical
line through the middle black oval. The other four points are obtained
from the first four by a rotation of order two around the line at height 
1/4 over any of the two horizontal axes of the big rectangle.

Let us also see how to check in the picture the conditions of the
statement. An isometry $\tau$ is in the normalizer of $G$ if and only
if for any generic orbit $GP$, $\tau GP$ is again an orbit of $G$ (if
this happens, let $P'=\tau P$ and let $P''$ be such that $\tau
GP=P''$. Then $\tau G\tau^{-1} P'=GP''$ which implies $P'\in GP''$ and,
since $P$ is generic, $\tau G\tau^{-1}=G$). In each figure it is easy
to check that a reflection in any of the displayed lines produces a
new orbit of $G$.

As for the number of planes in the band $Z_P$, each horizontal orbit
of $G$, except for the one containing $P$, produces exactly two planes
in $Z_P$, since $Z_P$ has width 2 and the minimal vertical translation
has length one. In the Figures, the left-top and right-bottom
quadrants are ``white'' fundamental subdomains and the other two are
``black''.  
\qed

\begin{figure}[htb]
{\centerline{
 \epsfxsize=9.5cm
 \epsfbox{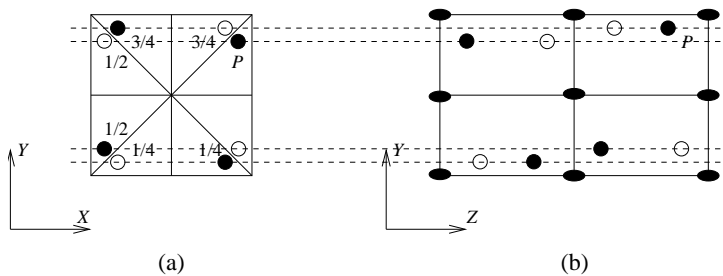}
}}
  \caption{The group $P4_122$}
\label{Fig.P4122}
\end{figure}

More or less the same technique can be applied to the group $P4_122$,
depicted in Figure \ref{Fig.P4122}.
Part (a) is its standard graphical representation,
showing the eight points of a generic orbit $GP$ which lie inside the
translational primitive cell of the tetragonal system, projected
onto the $XY$ plane and with their heights ($Z$-coordinates) indicated.
Part (b) shows the projection of these
same points to the $YZ$ coordinate plane. The black ovals
indicate order-2 rotations of $G$ with axes in lines
parallel to the $X$ axis. They
can also be read as order-2 rotation centers in the planar subgroup of $G$
which preserves planes parallel to $YZ$. This subgroup,
hence, is of type $p2$ with a rectangular grid.

\begin{Cor}
\label{Coro.P4122}
Dirichlet stereohedra for the group $P4_122$ cannot have more than 50 facets.
\end{Cor}

\Dem
Let $G$ be a group of type $P4_122$.
Let $G_1$ be the  subgroup of $G$ which preserves planes
parallel to $YZ$. Let $Z_P$ be the band centered at the $YZ$-plane
containing our base point $P$ and of width two times the shortest
translation of $G$ (and of $G_1$) in the $X$-direction. Only seven
$G_1$-orbits intersect the interior of the band. Although we cannot
say here that the $G_1$-orbits are related by the normalizer of $G_1$,
still we can apply Lemma \ref{Lema.ZP}, Lemma \ref{Lema.twoplanes} and 
Theorem \ref{Teo.dos.orbitas} to conclude that each $G_1$ orbit produces at
most 7 neighbors except $G_1P$ which produces at most 6. This gives at
most 48 neighbors in the interior of $Z_P$
and we have to add the two translates of $G_1$ in the
$X$-direction.

Observe that the choice of $P$ shown in Figure \ref{Fig.P4122}
produces four $G_1$-orbits in ``black'' fundamental subdomains, 
for which the number of
neighbors is at most 4. But there are choices of $P$ in which all the
$G_1$-orbits lie in white subdomains.
\qed

\subsection{Groups with a horizontal $p3$}

\begin{figure}[htb]
{\centerline{
 \epsfxsize=11cm
 \epsfbox{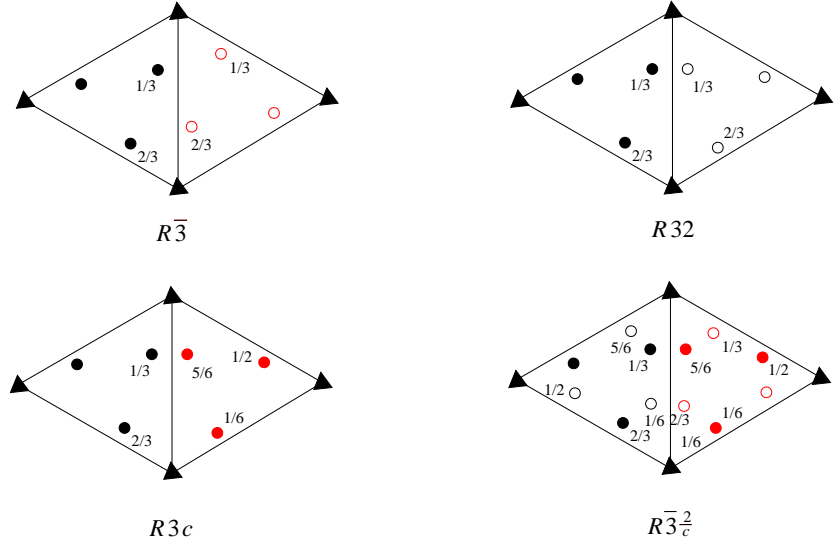}
}}
  \caption{Trigonal groups whose horizontal group is
of type $p3$}
\label{Fig.R3}
\end{figure}

\begin{Prop}
\label{Prop.R3}
The four trigonal groups displayed in Figure \ref{Fig.R3}
have a horizontal group of type $p3$. The first three have ten
horizontal planes (not counting the boundary and middle ones) in the region
$Z_P$, six of which lie over black subdomains and four over white
subdomains, in the sense of Lemma \ref{Lema.p3}.
Hence, their Dirichlet stereohedra cannot have more than 42 facets.

The last group has 22
horizontal planes (not counting the boundary and middle ones) in the region
$Z_P$, 12 which produce black orbits and ten which produce white
orbits. Hence, their Dirichlet stereohedra cannot have more than 84 facets.
\end{Prop}

\Dem As in Proposition \ref{Prop.rectangular.p2}, the stated bound follows
from the rest of the statement. The statement on horizontal planes follows
from the graphical representation of each group in Figure \ref{Fig.R3}. For
each group, the small black triangles represent centers  of rotation of order
three. The region displayed is a fundamental domain for the
horizontal group and the six (in the first three groups) or twelve (in the
last one) black or white points represent a generic orbit of $G$, with the
same conventions as in the proof of Proposition \ref{Prop.rectangular.p2}.
\qed

\begin{figure}[htb]
{\centerline{
 \epsfxsize=6 cm
 \epsfbox{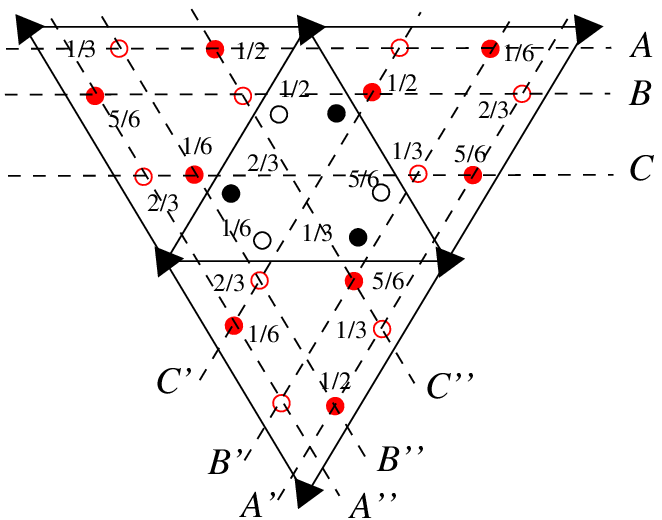}
}}
  \caption{Vertical planes in the group $R\overline{3}\frac{2}{c}$}
\label{Fig.R32c}
\end{figure}

For the trigonal group $R\overline{3}\frac{2}{c}$ the upper bound of 
84 neighbors obtained in Proposition \ref{Prop.R3} can still be
lowered a bit if we take into account that not only horizontal planes
but also vertical planes contain many points of each orbit.
Observe that the bound of Proposition
\ref{Prop.R3} can be stated more precisely as saying that $P$ can have
at most 48 neighbors which project to white triangles in the influence
region and at most 36 neighbors which project to black triangles. We
will be interested in the latter ones. Figure \ref{Fig.R32c} shows the
projection of a generic orbit to the union of 4 fundamental subdomains
forming a triangle. If we assume $P$ to project to the central
triangle, colored white, then the other three triangles are the three black
triangles contained in the influence region.

The dashed lines in Figure \ref{Fig.R32c} represent nine vertical planes
$A$, $B$, $C$, $A'$, $B'$, $C'$, $A''$, $B''$ and  $C''$ with the
property that any black neighbor of $P$ lies in two of them. On the
other hand, the subgroup of $G$ fixing each of those planes is a
planar crystallographic group of type $p2$ and, by Theorem
\ref{Teo.dos.orbitas}, produces at most seven neighbors. Hence the
number of black neighbors cannot exceed $9\times 7 / 2 =31.5$. This
implies the following: 

\begin{Cor}
\label{Coro.R32c}
A Dirichlet stereohedron for a group of type $R\overline{3}\frac{2}{c}$
cannot have more than 31+48=79 neighbors.
\end{Cor}

\section{Reduced influence region. Groups with a horizontal $pgg$ or
$pg$}
\label{sec.reduced}

\subsection{How two orbits of $pgg$ related by the normalizer intersect}
We recall the following result from \cite{MiArt1} (Theorem 3.1):

\begin{Lema}
\label{Lema.involutive.pgg}
Let $G_0$ be a planar crystallographic group of type $pgg$ and let $P$
and $Q$ be two points in the plane. Assuming that $G_0P\cup G_0Q$ is an orbit
of a certain crystallographic group, the number of Dirichlet regions
of $\Vor_{G_0Q}$ which overlap $\Vor_{G_0P}(P)$ is at most seven.
\noproof
\end{Lema}

A case study shows that under the hypotheses of the lemma we must have
$Q=\tau P$ where $\tau$ is an element of the normalizer
of $G_0$ in the group of Euclidean isometries of the plane. We intend
to extend the lemma to this more general case, thus proving part (vi)
of Theorem \ref{Teo.dos.orbitas}.

Recall that $pgg$ is generated by two perpendicular translations

together with a rotation of order two and a glide reflection of vector
half of one of the translations.
There are two possibilities for the normalizer:
\begin{itemize}
\item
If the generating translations have different lengths, then the
normalizer $N_0$ is generated by $G_0$ together with reflections on
the lines supporting the rectangles of the grid of rotation
centers. $G_0$ has index four in $N_0$ and we have that for any
$\tau\in N_0$, $\tau^2$ is in $G_0$. In particular, $G_0\cup \tau G_0$
is a crystallographic group for every $\tau\in N_0$, and hence
Lemma \ref{Lema.involutive.pgg} already implies what we
want to prove.
\item
If the generating translations have the same length, then the
normalizer $N_0$ is generated by $G_0$ together with the reflections
mentioned above and rotations of order four in the rotation centers of
$G_0$. This is the case we will be interested in.
\end{itemize}

\begin{figure}[htb]
  {\centerline{
        \epsfysize=5 cm
        \epsfbox{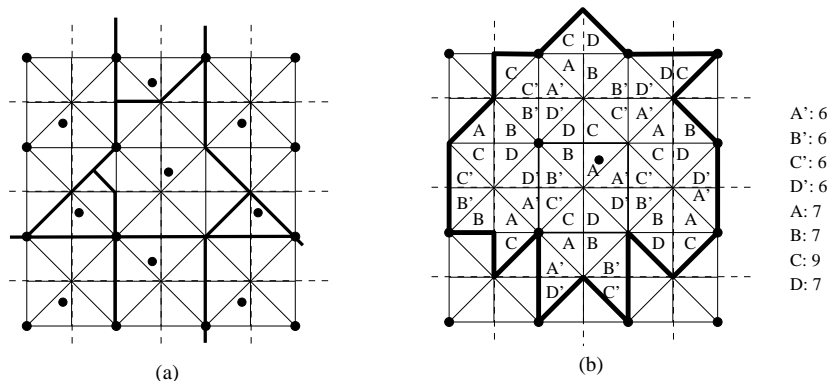}
}}
  \caption{Influence region of a $pgg$ with a square grid of rotation centers}
  \label{Fig.infl.pgg}
\end{figure}

$G_0$ has index eight in $N_0$. More precisely, we can take as fundamental
domains of $G_0$ squares with vertices in rotation centers and as
fundamental subdomains (i.e. fundamental domains of $N_0$) the
eight triangles in which the symmetries of the square
divide the fundamental domains. The influence region is computed in Figure
\ref{Fig.infl.pgg}. The left part is an extended Dirichlet region of
the initial fundamental subdomain $D$,
computed as in Figure \ref{Fig.p3}. The right part is the union of all
the fundamental subdomains whose extended Dirichlet region overlaps
the one on the left part. The extended Dirichlet region of a
fundamental subdomain $\tau D$ is obtained applying $\tau$ to the
extended Dirichlet region $\Ext_{G_0}(D)$, for each $\tau\in N_0$.

We have labeled the eight elements of $N_0/G_0$ (and hence the
fundamental subdomains) with the
letters A, B, C, D, A', B', C' and D'. A represents $G_0$ itself
and A, B, C and D would form the normalizer of $G_0$ in
case that the two generating translations had different length.
The number of
fundamental subdomains of each type 
in the influence region is shown
on the right part of Figure \ref{Fig.infl.pgg}. 
This number never exceeds nine. This implies:

\begin{Lema}
\label{Lema.pgg}
Let $G_0$ be a planar crystallographic
group of type $pgg$ whose generating translations have
equal length. Let $G_0P$ and $G_0Q$ be two orbits with trivial
stabilizer. Then, the number of Dirichlet regions of one of the orbits
overlapped by each Dirichlet region of the other orbit is at most nine.
\end{Lema}

\begin{figure}[htb]
  {\centerline{
        \epsfysize=5 cm
        \epsfbox{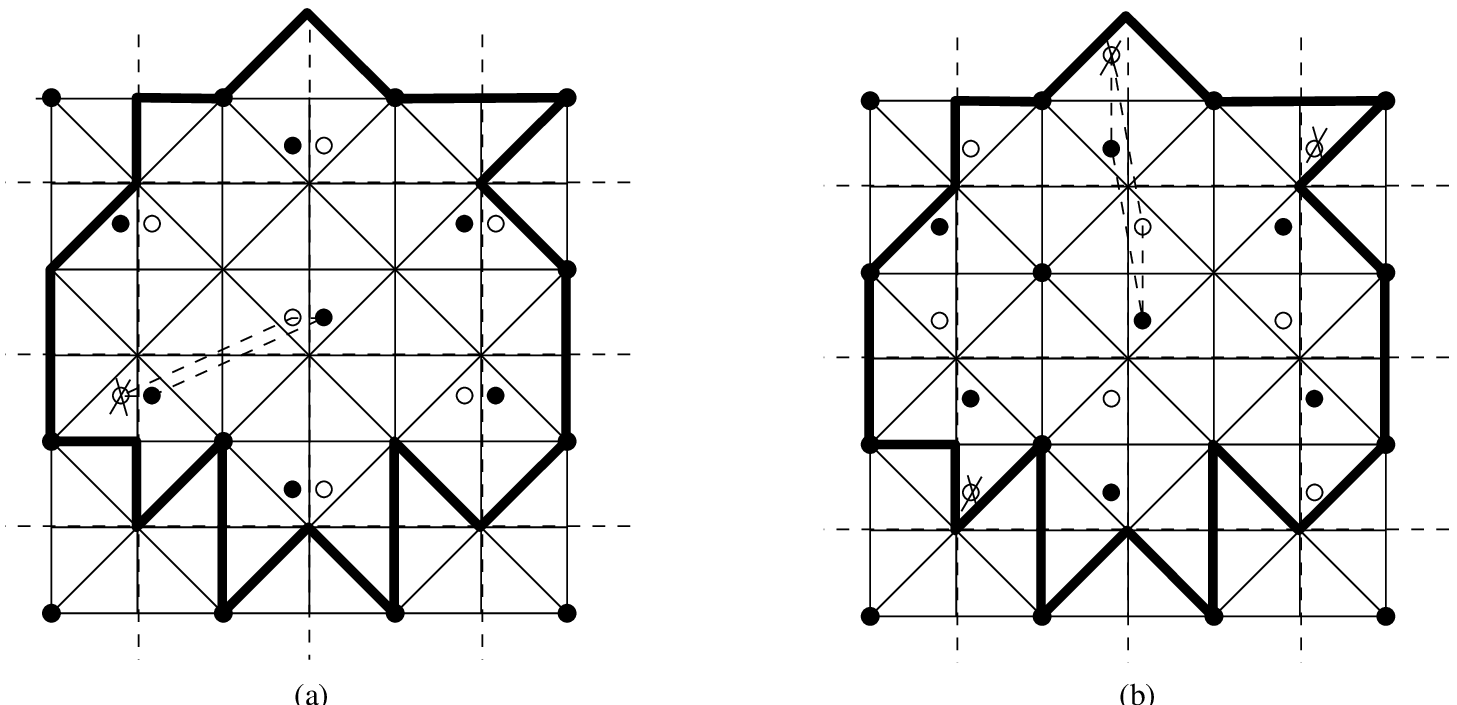}
}}
  {\centerline{
        \epsfysize=5 cm
        \epsfbox{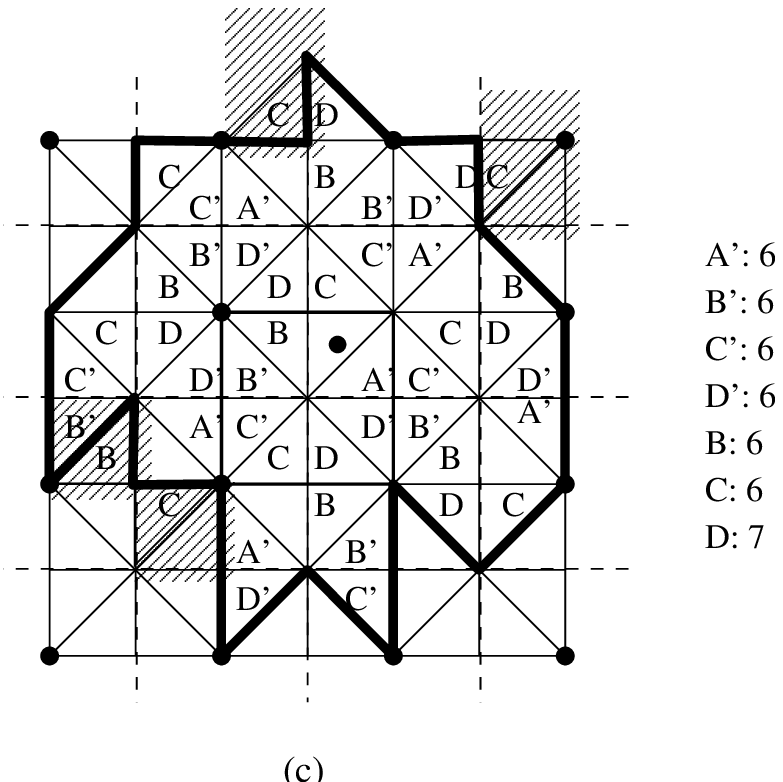}
}}
  \caption{Some fundamental subdomains in the influence region in
        Figure \protect{\ref{Fig.infl.pgg}} cannot
        produce overlapping Dirichlet regions, if the two orbits are
        related by the normalizer (parts (a) and (b)). This produces a
``reduced influence region'' (part (c)).}
  \label{Fig.red.pgg}
\end{figure}

So far we have not used the fact that we are interested in orbits
related by the normalizer. How to use this property is 
exhibited in Figure
\ref{Fig.red.pgg}. In part (a), our base orbit of $G_0$ is shown (in black)
together with the corresponding orbit of type ``B'' (in white). One point
of this latter orbit has been crossed out, meaning that its
corresponding Dirichlet region can never overlap the region of the
base point. The reason is in the four points joined by a dashed
quadrilateral: the movement of the normalizer sending the black
vertices of this quadrilateral to the white ones is an order-2 rotation, 
hence the quadrilateral is a parallelogram.
The bisectors of the two black vertices and the two white
vertices of the quadrilateral
will be parallel and will separate the Dirichlet regions of the
base point $P$ and the point crossed out.

The same argument applies to the point crossed out on the top end of
Figure \ref{Fig.red.pgg}(b), where the ``A'' and ``C'' orbits are shown.
But we have crossed out also two other points, one on the top-right end
and one on the bottom-left end. The proof that the Dirichlet regions
of these two points cannot overlap the one of our base point is in
\cite{MiArt1} (paragraph $[A,C]$ and, in particular, Figure 19).
Hence, instead of 7 and 9 fundamental subdomains of types B and C
which potentially could produce Dirichlet regions overlapping the base
Dirichlet region we have now 6 of each type. As a conclusion:

\begin{Prop}
\label{Prop.normal.pgg}
Let $G_0$ be a planar crystallographic group of type $pgg$ and let $P$
and $Q$ be any two points in the plane. Assuming that $Q=\tau P$ for some
element $\tau$ in the normalizer of $G_0$, the number of Dirichlet regions
of $\Vor_{G_0Q}$ which overlap $\Vor_{G_0P}(P)$ is at most seven.
\noproof
\end{Prop}

\subsection{Groups with a horizontal $pgg$}
The reduced influence region of Figure \ref{Fig.red.pgg} can be used to lower
the upper bounds given by Corollary \ref{Coro.planes} for the groups
whose horizontal subgroup $G_0$ is a $pgg$ with a square grid.
This follows from the fact
that only one type of planar orbits of those produced by
the normalizer give seven neighbors (the one we have labeled $D$),
and the rest only six.

More precisely, Figure \ref{Fig.tetra.pgg} shows the graphical
representation of the four groups in question.
The square displayed is a
fundamental domain of $G_0$, divided into eight fundamental subdomains.
In the first three groups the ``bad''
planar orbit, labeled with a D in Figure
\ref{Fig.infl.pgg}, does not appear. Hence we can take $i=6$ in the
computations of Corollary \ref{Coro.planes} for these groups.
In $I\frac{4_1}{g}\frac{2}{c}\frac{2}{d}$ the bad orbit appears, but
still we can count six neighbors for 12 of the planes in the band
$Z_P$ and seven for only two of them. This gives:

\begin{figure}[bht]
{\centerline{
 \epsfxsize=11cm
 \epsfbox{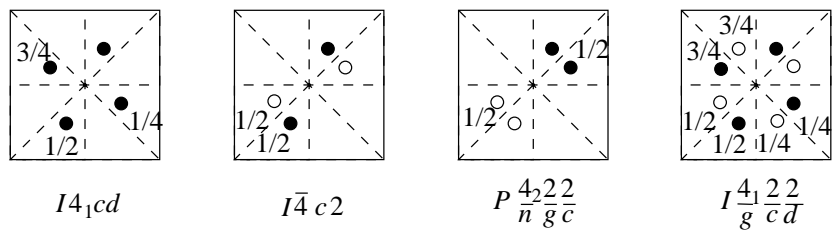}
}}
  \caption{Tetragonal groups whose horizontal group is
of type $pgg$}
  \label{Fig.tetra.pgg}
\end{figure}

\begin{Prop}
\label{Prop.pgg}
Dirichlet stereohedra for the groups $I4_1cd$, $I\overline{4}c2$ and
$P\frac{4_2}{n}\frac{2}{g}\frac{2}{c}$ cannot have more than 44
facets.
Dirichlet stereohedra for the group
$I\frac{4_1}{g}\frac{2}{c}\frac{2}{d}$ cannot have more than 94
facets.
\end{Prop}

\begin{figure}[bht]
{\centerline{
 \epsfxsize=12cm
 \epsfbox{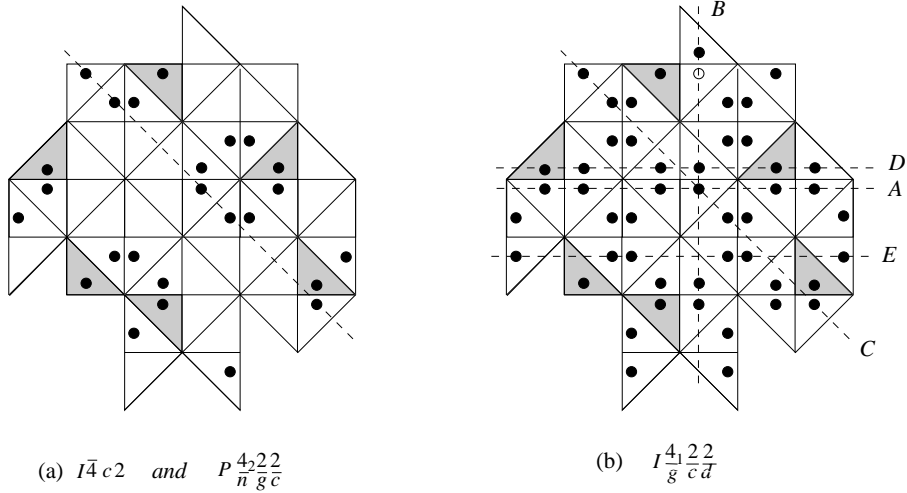}
}}
  \caption{Vertical planes in the groups 
of type $pgg$}
  \label{Fig.pgg.vert}
\end{figure}

But we can use vertical planes to refine this result a bit.
Figure \ref{Fig.pgg.vert} shows the projection of
a generic orbit $GP$ to the influence region computed above, where $G$
is of type $I\overline{4}c2$ or $P\frac{4_2}{n}\frac{2}{g}\frac{2}{c}$
in part (a) of the Figure and of type
$I\frac{4_1}{g}\frac{2}{c}\frac{2}{d}$ in part (b). The number
44 or 94 
stated in Proposition \ref{Prop.pgg} is obtained counting 2 possible
neighbors over each fundamental subdomain, except in
the six grey subdomains (those containing points at the same height
as $P$ except the one containing $P$ itself) where only one neighbor
over each fundamental subdomain is possible. 
The dashed line in part (a) of the figure represents a vertical plane
containing the base point $P$ where we have counted ten
possible neighbors while there are at most six real neighbors (since
in a planar Dirichlet tiling each region is a neighbor of at most
other six). Hence, instead of 44 we can take 40 as an upper bound
for $I\overline{4}c2$ or $P\frac{4_2}{n}\frac{2}{g}\frac{2}{c}$.

Similarly, the five dashed lines 
in part (b) of Figure \ref{Fig.pgg.vert}
represent five vertical planes $A$, $B$, $C$, $D$ and $E$ 
where we have counted 12, 12, 10, 10 and 12 possible
neighbors respectively, giving a total of 50
because two of them are common to $A$, $B$ and $C$, and other two
common to $B$ and $D$.
As before, $A$, $B$ and $C$ can produce at most six neighbors each,
because they all contain $P$.

Let $G_1$ denote the subgroup of $G$ preserving planes parallel to
$A$, in particular preserving $D$ and $E$. $G_1$
is of type $pgg$ (although it does not have its two translations
of equal length). According to Theorem \ref{Teo.dos.orbitas} each
plane parallel to $A$ can contain at most 11 neighbors, and only seven
if the $G_1$-orbits in that plane and in $A$ are related by the
normalizer. The latter happens in the plane $D$; the orbits in the
planes $A$ and $D$ 
are related by a glide reflection in the bisecting plane, which
projected to any of them becomes a translation lying in the normalizer of
the corresponding $pgg$.

Hence, the planes $A$, $B$, $C$, $D$ and $E$ can contain in total at
most $6+6+6+11+7=36$ neighbors, instead of the 50 which we have counted.

\begin{Cor}
\label{Coro.pgg.vert}
Dirichlet stereohedra for the groups 
$I\overline{4}c2$ and $P\frac{4_2}{n}\frac{2}{g}\frac{2}{c}$
cannot have more than 40 facets. Those for 
$I\frac{4_1}{g}\frac{2}{c}\frac{2}{d}$ cannot have more than 80 facets.
\end{Cor}

\begin{Obs}
\rm
\label{78}
Using the methods of \cite[Section 3.3]{MiArt1} in a more sophisticated
way, the number appearing for $pgg$ in part (v) of Theorem
\ref{Teo.dos.orbitas} can be lowered to nine instead of eleven
(see \cite[Remark 3.2]{MiArt1}). In the
preceeding argument, this would
lower the number of possible neighbors in the plane $E$ 
by two as well, giving a bound of 78 instead of 80. But this small
improvement does not seem to
be worth the effort of giving a proof here.
\end{Obs}

\subsection{Groups with a horizontal $pg$}

\begin{figure}[htb]
  \centerline{
    \epsfysize=8cm
    \epsfbox{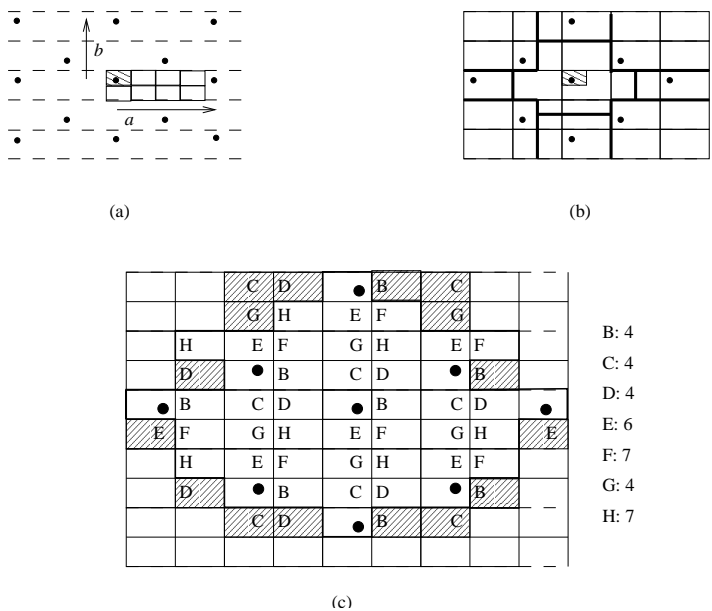}
}
\caption{
    Extended Dirichlet region and influence region for the group $pg$}
    \label{Fig.pg}
\end{figure}

We will now compute the influence region of a planar group of type
$pg$, in order to lower the bound for the groups
$P\frac{2_1}{n}\frac{2}{c}\frac{2_1}{a}$ and
$P\frac{2_1}{a}\frac{2_1}{b}\frac{2_1}{c}$.
Remember that $pg$ is
generated by two perpendicular translations and a glide reflection
on a line parallel to one of them. Let $a$ and $b$ be the vectors of
the two translations, $a$ being parallel to the glide-reflection
line. Observe that $G$ consists of translations and glide-reflections
on a family of lines parallel to $a$ and distant $\frac{|b|}{2}$ to one
another. Any rectangle of sides $a$ and $\frac{b}{2}$ placed between two
consecutive such lines is a fundamental domain for $G_0$. See Figure
\ref{Fig.pg}(a).

We take as $N_0$ the group generated by $G_0$, a reflection on a line
parallel to $b$ and reflections on the glide-reflection lines of $G_0$
and the midlines between any two consecutive glide reflection lines.
Each fundamental domain of $G_0$ gets divided
into 8 fundamental domains of $N_0$, which are rectangles of sides
$\frac{a}{4}$ and $\frac{b}{4}$. See again Figure \ref{Fig.pg}(a), where
one of these fundamental subdomains has been shaded.

Part (b) of Figure \ref{Fig.pg} shows the extended Dirichlet region of
a fundamental subdomain, and part (c) the corresponding influence
region. The shadowed rectangles in part (c) are the fundamental
subdomains lying in the influence region but which cannot produce
overlapping regions if the two orbits of $G_0$ are related by $N_0$,
with the same argument as in part (a) of Figure
\ref{Fig.red.pgg}. Hence, the interior of the thick polygon is the
reduced influence region. We label the eight $G_0$-cosets
in $N_0$ with the letters A, B, C, D, E, F, G and H where A is $G_0$
itself, and get for the orbits in each class the number of overlapping
regions shown on the right of Figure \ref{Fig.pg}(c). This number
equals four in the cosets obtained from A by a translation in the
direction of $a$ (coset C), reflection on a line parallel to $b$
(cosets B and D), or glide-reflection with axis in the direction of $a$
(coset G). It equals six in the coset E
obtained by any of the other of reflections of $N_0$, and seven in
the other two cosets F and G, obtained by order-two rotations
or vertical translation.

\begin{figure}[htb]
{\centerline{
 \epsfysize=2cm
\epsfbox{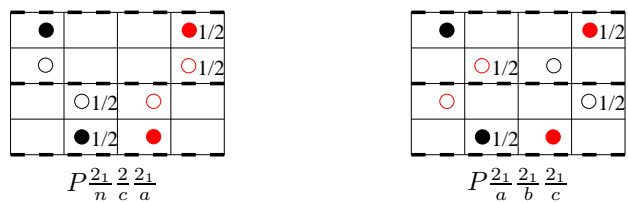}
}}
{\centerline{
$P\frac{2_1}{n}\frac{2}{c}\frac{2_1}{a}$ \hspace{3.87cm}
$P\frac{2_1}{a}\frac{2_1}{b}\frac{2_1}{c}$
}}
\caption{Orthorhombic groups with the horizontal planar group $pg$}
\label{Fig.pg.ortho}
\end{figure}

We now look at the two orthorhombic groups in Table \ref{Table.Orthorhombic}
having a horizontal $pg$. They are depicted in Figure
\ref{Fig.pg.ortho}. As usual, the picture shows the projection of a
translational cell, which in this case consists of two fundamental
domains of $G_0$. In the group $P\frac{2_1}{n}\frac{2}{c}\frac{2_1}{a}$
we see that the horizontal planes contain orbits of types D, E and H,
while in $P\frac{2_1}{a}\frac{2_1}{b}\frac{2_1}{c}$ we have of types
D, F and G. 

\begin{Prop}
\label{Prop.pg}
Dirichlet stereohedra for the groups
$P\frac{2_1}{n}\frac{2}{c}\frac{2_1}{a}$  and
$P\frac{2_1}{a}\frac{2_1}{b}\frac{2_1}{c}$ cannot have more than 
38 facets, respectively.
\end{Prop}

\Dem 
In $P\frac{2_1}{a}\frac{2_1}{b}\frac{2_1}{c}$, the band $Z_P$ contains
two horizontal orbits of each of types D, F and G, which produce at
most $2\times(4+7+4)=30$ neighbors. These, added to the six neighbors
in the base horizontal orbit and the two vertical translates of $P$,
provides the upper bound.

\begin{figure}[htb]
  \centerline{
    \epsfysize=5cm
    \epsfbox{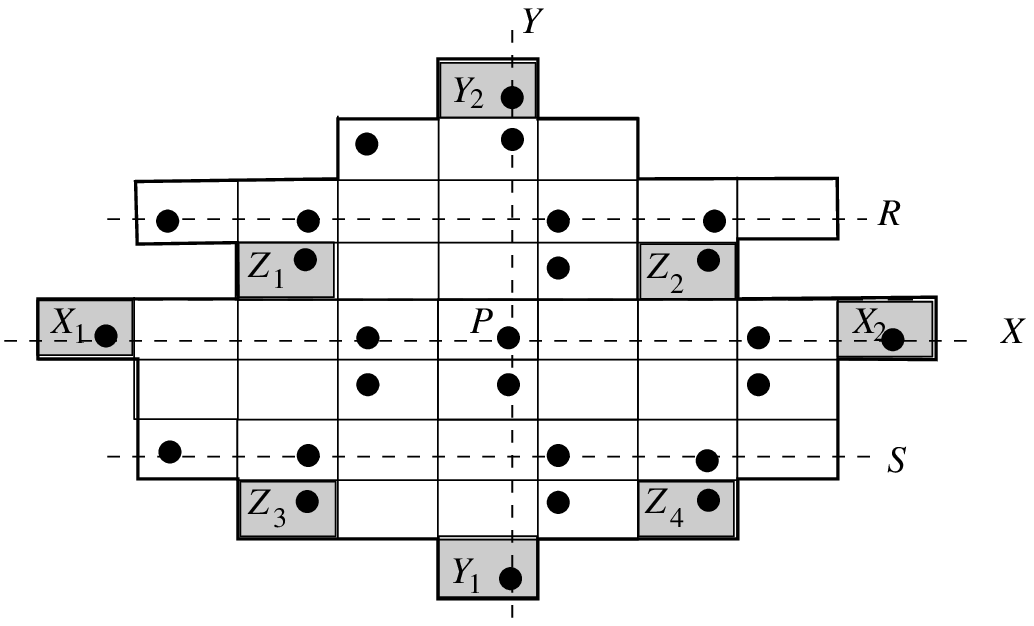}
}
\caption{
    Vertical planes in the group $P\frac{2_1}{n}\frac{2}{c}\frac{2_1}{a}$}
    \label{Fig.pg.vert}
\end{figure}

For the group
$P\frac{2_1}{a}\frac{2_1}{b}\frac{2_1}{c}$ we will use vertical
planes, and
Figure \ref{Fig.pg.vert} which shows the projection of a generic orbit of
$G$ to the reduced influence region computed above. 

Counting 2 points over each white rectangle
and one over each grey rectangle gives 44 possible neighbors. 

Let us concentrate in the $G_0$-orbit containing $P$, whose
intersection with the reduced 
influence region consists of $P$ and the eight
points in grey rectangles, labeled $X_1$, $X_2$, $Y_1$,
$Y_2$, $Z_1$, $Z_2$, $Z_3$ and $Z_4$ in Figure \ref{Fig.pg.vert}.
$X_1$, $X_2$, $Y_1$ and
$Y_2$ are obtained from $P$ by translations of $G_0$. The other four
by glide reflections.
In the planar Dirichlet tiling produced by this orbit we have:

\begin{enumerate}
\item[(i)] $P$ is a neighbor of the point $X_1$ if and only if it is a
neighbor of $X_2$, because $P$
being a neighbor of $Q$ implies $\tau P$ being a neighbor of $\tau Q$
for any $\tau\in G_0$.

\item[(ii)] $P$ is a neighbor of the point $Y_1$ if and only if it is a
neighbor of $Y_2$, for the same reason.

\item[(iii)] The four points $X_1$, $X_2$, $Y_1$ and $Y_2$ cannot be
all neighbors of $P$. This is a topological argument: If $X_1$ and
$X_2$ are neighbors of $P$ then $Z_1$ is a neighbor of $Z_2$ and 
$Z_3$ is a neighbor of $Z_4$. There
is no ``room left'' for $Y_1$ or $Y_2$ being neighbors of $P$.
\end{enumerate}

Let us say that we are in the ``X'' case if $Y_1$ and $Y_2$ are
not neighbors of $P$ and in the ``Y'' case if $X_1$ and $X_2$ are
not neighbors of $P$. Statement (iii) above implies that we are always
in one of the two cases (maybe in the two of them).

The dashed lines in Figure \ref{Fig.pg.vert} represent
four vertical planes $R$, $S$, $X$ and $Y$. The number
44 implied counting 8 possible neighbors in each of them and 30 points
in total in the four planes, since the two points counted over $P$ are
common to $X$ and $Y$. We will see that the four planes can contain
at most 26 neighbors of $P$ in total, which finishes the proof.

$R$ and $S$ can contain at most seven neighbors
each by Theorem
\ref{Teo.dos.orbitas}, since the planar subgroup
of $G$ with respect to the planes parallel to $R$, $S$ and $X$ is of
type $pg$.

$X$ and $Y$ contain $P$ and hence each of them gives at most six
neighbors. But we can be more precise: in the ``X'' case, $X$ gives
at most six neighbors and $Y$ gives at most four neighbors not
contained in $X$. In the ``Y'' case, $Y$ gives
at most six neighbors and $X$ gives at most four neighbors not
contained in $Y$.
\qed

\section{The group $P6_122$}

\subsection{An upper bound using vertical planes}

\begin{figure}[htb]
  \centerline{
        \epsfxsize=10 cm
        \epsfbox{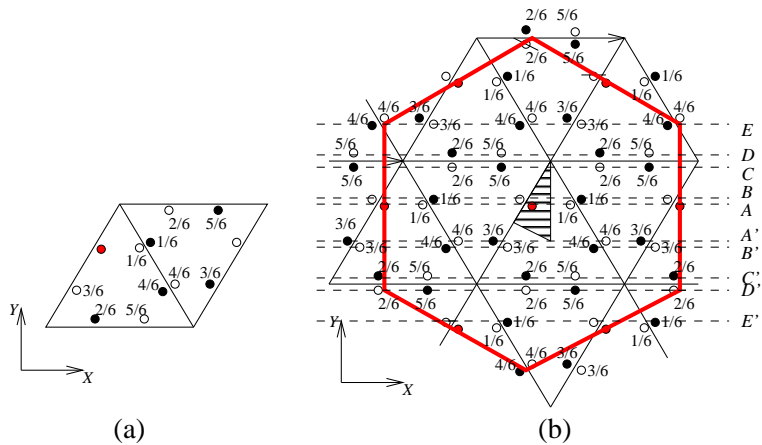}
  }
  \caption{$P 6_1 22 $:  (a) Neighbors in the interior of $Z_P$ 
          (b)  $XZ$ view }
\label{Fig.p1.P6122}
\end{figure}

We now deal with a group $G$ of type $P6_122$, whose horizontal
subgroup $G_0$ is generated by two translations of equal length
forming an angle of $60^\circ$. The whole group is generated by $G_0$
together with 
\begin{itemize}
\item a screw rotation of order 6 and of vertical axis and
\item any rotation of order 2 with axis parallel to one of the
generating horizontal translations and intersecting the screw-rotation
axis.
\end{itemize}
See a graphical
representation of the group in Figure \ref{Fig.p1.P6122}(a).

\begin{Prop}
\label{Prop.P6122}
Dirichlet stereohedra for the group
$P6_122$ cannot have more than 78 facets.
\end{Prop}

\Dem
Remember that Corollary \ref{Coro.planes} gave an upper bound of 96
for the number of neighbors in this group. We first show a
different way of deriving this same upper bound which will be more
appropriate for our purposes here.

Let $l$ denote the minimal length of a horizontal translation in $G$, so that
the horizontal group $G_0$ generates a triangular lattice with equilateral
triangles of side $l$. Each such triangle has height $a=\sqrt{3}l/2$. 
Then, for any point
$P$ we have that $\Vor_{G_0P}(P)$ is an infinite prism over a regular hexagon 
of side $2a/3=l/\sqrt{3}$. Since this happens for every point,
the condition necessary and sufficient for 
$\Vor_{G_0P}(P)\cap \Vor_{G_0Q}(Q)\ne\emptyset$ is than that $Q$ lies
in the prism over $2\Vor_{G_0P}(P)$, i.e., over the regular hexagon of
side $4a/3=2l/\sqrt{3}$ centered at $P$ and with sides orthogonal to
the primitive translations of $G_0$. This is the thick hexagon
in Figure \ref{Fig.p1.P6122}(b). Moreover, if $Q$ is not in
the same horizontal plane as $P$, a necessary condition for a $Q\in GP$ being
a neighbor of $P$ in $\Vor_{GP}(P)$ is that the projection of $Q$ to
the horizontal plane containing $P$ lies strictly inside that
hexagon. For the points in the same horizontal plane as $P$ it is
allowed, however, to lie in the boundary of the hexagon. Then, the
possible neighbors of $P$ in $\Vor_{GP}(P)$ are:

\begin{itemize}

\item The six points obtained from $P$ by primitive horizontal
translations, i.e., the mid-points of the edges of the thick hexagon in Figure
\ref{Fig.p1.P6122}.

\item The two closer vertical translates of $P$.

\item The points which project to the interior of the thick hexagon
and with vertical distance to $P$ smaller than 1. There are 22
horizontal planes other than the one containing $P$ and with distance
to $P$ smaller than 1. A simple density argument shows that
at most four orbit points in each of the 22 planes lie inside the hexagon. More
precisely, the density argument says that {\em exactly} four points at
each horizontal plane project to the hexagon if one counts points
projecting to the interior as 1, those projecting to
facets of the hexagon as $1/2$ and those projecting to
vertices as $1/3$.
\end{itemize}

Summing up, the above gives a bound of $22\times4+6+2=96$
neighbors. Our goal is to decrease the $22\times 4$ part by
considering vertical planes parallel to the $XZ$ coordinate plane.
The subgroup $G_1$ consisting of elements of $G$ which preserve
those planes is of type $p2$, generated by a translation of length $l$
in direction $X$, a translation of length 1 in direction $Z$ and any
order 2 rotation contained in $G_1$ (whose axis will be parallel to
$Y$ and intersect some screw rotation axes of $G$).

The orbit $GP$ decomposes into an infinite family of orbits of the group
$G_1$ lying in different $XZ$-planes. Ten of them are marked in
Figure \ref{Fig.p1.P6122}(b) with the letters $A$, $B$, $C$, $D$, $E$,
$A'$, $B'$, $C'$, $D'$ and $E'$. $A$ is the one containing the
base point $P$. In each pair $A$-$A'$, $B$-$B'$, $C$-$C'$, $D$-$D'$
and $E$-$E'$ one of the $G_1$ orbits is obtained from the other by a
screw rotation $\rho$ of order 2 in the vertical axis $\{X=-l/4,
Y=-a/2\}$. More importantly, when projected to any of the $XZ$-planes,
$\rho$ becomes an element of the normalizer of $G_1$ which exchanges
the coloring of the rectangular tiling used  in Lemma
\ref{Lema.rectangular.p2}. Hence, by that Lemma, 

\bigskip
{\bf Claim 1:} In $A'$ we can have at
most 4 neighbors of $P$ and in any of the other pairs $B$-$B'$, $C$-$C'$,
$D$-$D'$ and $E$-$E'$ one of the two planes can provide at most 7
neighbors and the other one at most 4 neighbors (although a priori we
do not know which plane is which). 
\bigskip

\begin{figure}[htb]
  \centerline{
        \epsfxsize=6 cm
        \epsfbox{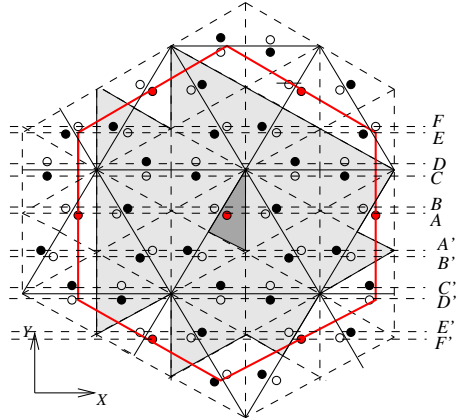}
  }
  \caption{Vertical planes in the group $P 6_1 22 $}
\label{Fig.p1.P6122.2}
\end{figure}

Let us now show that the bound of 96 obtained above over-counted the
points in these ten planes by at least 18. For this we use
Figure \ref{Fig.p1.P6122.2}, which essentially coincides with Figure
\ref{Fig.p1.P6122}(b) except that we have removed the numbers showing
the height of different points and included instead the subdivision of
the plane into fundamental subdomains of the horizontal group
$G_0$. The darker fundamental subdomain, which we call $D_0$, is the
one containing the base point $P$.

The lighter shaded
polygonal region of Figure \ref{Fig.p1.P6122.2} consists of the
fundamental subdomains $D$ with the property that for any choice of $P$ in
$D_0$ the points of $GP$ lying over $D$ are in the interior of the
hexagonal influence region centered at $P$. In other words, the
fundamental subdomains $D$ for which the upper bound of 96 includes two
points above $D$ for any choice of $P$.

This implies that the 96 points of the previous bound include at least
6 points in the planes $D'$, $C'$ and $E$, at least 
8 points in the planes $B'$, $A'$, $A$, $B$, $C$ and $D$ and at least
5 points in the plane $E'$. The number for $A$ include the special
counting of points obtained from $P$ by translations of $G$. The
number for $E'$ takes into account a fundamental subdomain whose orbit
points lie in a facet of the hexagonal prism and which have
been counted as 1/2 in the ``96 points''. According to
Claim 1, we have over-counted at least 5, 3 and 3,
respectively, in the pairs of planes $B$-$B'$, $C$-$C'$ and $D$-$D'$.
In the pair $E$-$E'$ we have over-counted at least 1, since one of them
produces at most 4 neighbors and we have counted more than 4 in each
of them.
Finally, we have over-counted 4 points in the plane $A'$ and 2 in
the plane $A$.
\qed

\subsection{Dirichlet stereohedra with many facets}

Our constructions of Dirichlet stereohedra with many facets are based on
the following result. The case $k=8$ appeared in \cite{MiArt1}. The general
case has been independently proved by J.~Erickson as Lemma 1 in 
\cite{Erickson}:

\begin{Lema}
\label{Lema.screw}
Let $g$ be the screw rotation of order $k$ around the $Z$ coordinate axis, 
with translation of length l, i.e.
$g(x,y,z)=(x \cos(2\pi/k) - y \sin(2\pi/k), x \sin(2\pi/k) + y \cos(2\pi/k),
z+l)$.
Let $\langle g\rangle$ be the (infinite cyclic)
group generated by $g$. Then, any point $P$ not on the
axis of $g$ has as neighbors in $\Vor_{\langle g\rangle P}(P)$ the
$2k$ points $g^i P$ for $i \in\{-k,\dots,k\}\setminus \{0\}$.
\end{Lema}

\Dem
For any $P\in \R^3$ (but not on the $Z$-axis)
the convex hull of $\langle g\rangle P$ is an infinite prism
over a regular $k$-gon, one of the edges of the prism containing the
points $g^{nk} P$, for $n\in\Z$. In particular, $g^k P$ and $g^{-k}(P)$
must be neighbors of $P$. We will prove that the points $g^i P$ for
$i\in \{ -k+1,\dots,k-1\}\setminus\{0\}$ are also neighbors of $P$.

For any $P\in \R^3$ the orbit $\langle g\rangle P$ is contained in
an helicoidal curve $H_{r,\alpha,l}:=\{(r\cos(t-\alpha),
r\sin(t-\alpha), \frac{l}{2\pi}t)\in\R^3 : t\in \R\}$. There is no loss of
generality in fixing $l=2\pi$ and $\alpha=0$, hence using the
notation $H_r$ for $H_{r,0,2\pi}$.
We denote $P_t=(r\cos(t), r\sin(t), t)$ the point of $H_r$
for a certain parameter $t$, so that $g^i P_t = P_{t+2\pi i/k}$. 
We will prove that
for any two points $P_t$ and $P_{t'}$ whose angular distance $|t-t'|$
is less than $2\pi$, there is a sphere tangent to $H_{r}$ at these two
points and not containing any other point of $H_{r}$. This
implies the lemma.

To prove our claim, 
by symmetry considerations we can further assume that $t'=-t$. 
Then, the two points in the claim are of the form 
$P_{t_0}=(r\cos(t_0),r\sin(t_0), t_0)$ and 
$P_{-t_0}=(r\cos(t_0), -r\sin(t_0), -t_0)$ for a certain $t_0\in (0,\pi)$.
The sphere tangent to $H_{r}$ at $P_{t_0}$ and $P_{-t_0}$ will have
center at a point $O=(a,0,0)$ which satisfies that the vectors 
\[
OP_{t_0}=(r\cos(t_0)-a,r\sin(t_0), t_0)
\]
and 
\[
(dP_{t}/dt)_{t=t_0}=(r\sin(t_0), -r\cos(t_0), 1)
\]
are orthogonal.
The equation
\[ 
(r\cos(t_0)-a)r\sin(t_0) - r^2\cos(t_0)\sin(t_0) + t_0 =0
\]
gives the solution
\[
a=\frac{-t_0}{r\sin(t_0)}
\]

We now prove that for any point $P_t$ other than $P_{t_0}$ or
$P_{-t_0}$ the distance from $O$ to $P_t$ is strictly bigger than the
distance from $O$ to $P_{t_0}$. It suffices to consider
$t\in [-\pi,\pi]$ 
because $d(O,P_{t \pm 2\pi}) > d(O,P_t)$ for any
$t\in[-\pi,\pi]$. Moreover, 
 since $d(O,P_t) = d(O,P_{-t})$ we restrict our attention
to $t\in [0,\pi]$.
We define the function
\[
f(t):=d(O,P_t)^2-d(O,P_{t_0})^2 =\frac{2 t_0
       (\cos(t)-\cos(t_0))}{\sin(t_0)} +
       t^2-{t_0}^2
\]
whose first and second derivatives are
\[
f'(t)=2t - 2\frac{t_0}{\sin(t_0)}\sin(t),\qquad
f''(t)=2 - 2\frac{t_0}{\sin(t_0)}\cos(t).
\]

We have 
\begin{itemize}
\item[(i)] $f(t_0)=f'(t_0)=0$ (as expected).
\item[(ii)] $f''(t_0)>0$. For this, observe that $f''(t_0)>0$ is equivalent
to $\sin(t_0)-t_0\cos(t_0)>0$, which holds because the function 
$g(t):= \sin(t)-t\cos(t)$ is zero at the origin and its derivative
$g'(t)= t\sin(t)$ is strictly positive in $(0,\pi)$.
\item[(iii)] $f'$ is injective in $(0,\pi)$. 
Indeed, if
$f'(t_1)=f'(t_2)$ then $\sin(t_1)/t_1= \sin(t_2)/t_2$ and, in
particular, the derivative of the function $h(t):=\sin(t)/t$ must have
a zero between $t_1$ and $t_2$. But $h'(t)=(t\cos(t)-\sin(t))/t^2=0$
would imply $t=\tan(t)$, which does not happen in the interval $(0,\pi)$.
\end{itemize}

With this we prove that $f$ achieves its unique minimum in the
interval $[0,\pi]$ at $t=t_0$, as follows: 
Claim (iii) implies that $f''$ is either always non-negative or
always non-positive. By (ii) it is always
non-negative. Claims (i) and (ii) imply that $f$ has a local minimum at
$t_0$ which, by the previous observation, is the unique global minimum.
  \qed

Let $G$ be a crystallographic group of type $P6_122$. 
Its two metric
parameters (which define $G$ modulo conjugation by an isometry) are
the lengths $l$ and $a$ of its minimal translations in directions
parallel and perpendicular, respectively, to the order-6
screw-rotation axes. As usual, we assume that the
screw rotation axes are vertical, that is to say, parallel to the third
coordinate axis. Also, assume that one of the minimal horizontal
translations is parallel to the $X$-axis and that the $X$-axis itself
is one of the horizontal order-2 rotation axes in $G$. This completely
specifies the group $G$ (and agrees with Figures \ref{Fig.p1.P6122}
and \ref{Fig.p1.P6122.2}).

Let $g$ be the screw rotation of order twelve obtained substituting
$k=12$ in Lemma \ref{Lema.screw}. Observe that $g^2\in G$.
Let $\rho$ be the order-2 rotation
$(x,y,z)\mapsto (x,-y,-z)$ which is in $G$ by hypothesis. 
Suppose now that our base point $P$ has coordinates $(r \cos(\pi/12),
r\sin(\pi/12), l/24)$, so that $\rho(P)=g^{-1}(P)$. Under these
assumptions we have that the orbits of $P$ under the groups 
$\langle g\rangle$ and $\langle \rho,g^2\rangle \subset G$
coincide. By Lemma \ref{Lema.screw}, the Voronoi diagram of this part
of $GP$ alone produces $24$ neighbors of $P$. Now, if we fix the
parameters $l$ and $r$ in the above description and make $a$ tend to
infinity, the neighbors of $P$ in the Voronoi diagram of 
 $\langle \rho,g^2\rangle P$ will keep being neighbors in
$GP$. (Observe that $GP$ is obtained as the Minkowski sum of
$\langle \rho,g^2\rangle P$ and a triangular grid of side $a$ in a
horizontal plane). Moreover, since the regions in the Voronoi diagram
of $\langle \rho,g^2\rangle P$ are unbounded, new neighbors are
guaranteed to appear. Hence, the above construction is guaranteed to
produce at least 25 neighbors.

We could try to argue geometrically how many ``new'' neighbors have to
appear in the limit of $a$ going to infinity. For example, it is
relatively easy to show that this number does not depend on the choice
of $l$ and $r$. 
However, it seems easier to compute that number
experimentally:

\begin{Exm}
\label{Exm.31}
Taking $l=12$, $r=1$ and $a=100$ in the above setting
(which gives $P=(\cos(\pi/12),\sin(\pi/12), 1/2)$) one gets 31
neighbors. The number and identity of the neighbors is
(experimentally) stable under increasing the value of $a$.
\end{Exm}
Even more, playing with different possibilities for the
parameters we have found that:

\begin{Exm}
\label{Exm.32}
The metric parameters $l=100$, $a=950$ and the base point 
$P=(1, \tan(\pi/12), 4)$ produce a Dirichlet stereohedron with 32
neighbors.
This stereohedron is very unstable. For
example, changing the last coordinate of $P$ to be $100/24=4.166$,
which would match exactly the above setting with $r= 1/\cos(\pi/12)$, only
30 neighbors are obtained (and their identity changes drastically).
Also, the
parameter $a$ is not big enough for the number of neighbors to be
the same for bigger values of $a$.
\end{Exm}

\begin{figure}[htb]
  \centerline{
        \epsfxsize=10 cm
        \epsfbox{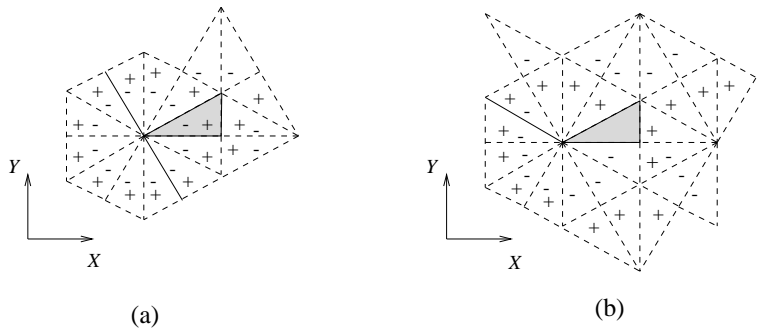}
  }
  \caption{Neighbors of $P$ obtained in Examples \ref{Exm.31} and
        \ref{Exm.32} (parts (a) and
        (b) respectively, for a group $P6_122$.}
\label{Fig.manyfacets}
\end{figure}

Figure \ref{Fig.manyfacets} describes the points of the orbit which
produce neighbors in each of the two examples. 
It shows the same decomposition of the
plane into fundamental subdomains which appears in Figure
\ref{Fig.p1.P6122.2}. The grey subdomain is the one containing $P$.
The points of the 
helicoidal orbit $\langle \rho,g^2\rangle P$ are those projecting to
the regular hexagon.
A plus (resp. a minus) in a subdomain means that the point of $GP$
projecting to that subdomain and lying in
the upper (resp. lower) half of the band $Z_P$ 
is a neighbor of $P$. In the first example all the 24
possible neighbors within the hexagon appear, as predicted. In the
second example only 16 of them appear, but this is compensated by 
more neighbors out of the hexagon.

The computations were made using a Maple program which first generates
a sufficiently large number of points $S$ of the orbit $GP$ to
guarantee that all the neighbors of $P$ lie in $S$ and then checks for
each of them whether it is actually a neighbor or not. For the first
step, it would be enough for example to use the two points in $Z_P$ in
each of the subdomains which appear in Figure
\ref{Fig.p1.P6122.2}. For the second step, we express being a
neighbor of $P$ as feasibility of a certain linear program.

The same ideas can be applied to groups
having an order-4 screw rotation around a
vertical axis and an order-2 rotation in a horizontal
axis, with the two axes intersecting one another.
One has to use Lemma \ref{Lema.screw} with $k=8$ and, 
hence, can only guarantee
to obtain more than 16 neighbors. This is the way
we constructed stereohedra with 18 facets for the group $I \frac{4_1}{g}
\frac{2}{m}\frac{2}{d}$, hence showing that this is exactly the highest
possible number of facets of Dirichlet stereohedra for groups with
reflections (see \cite[Example 2.9]{MiArt1}).

\begin{figure}[htb]
  \centerline{
        \epsfxsize=10 cm
        \epsfbox{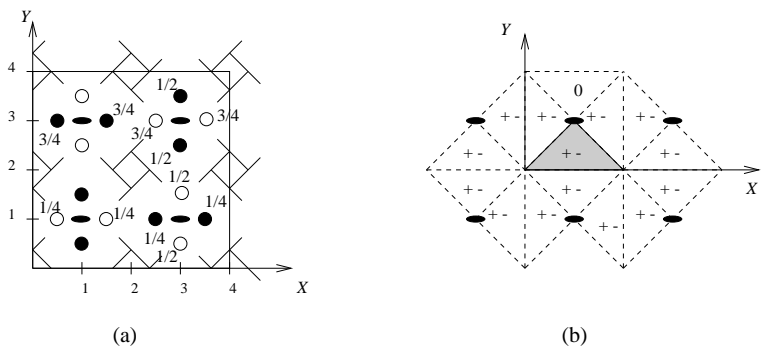}
  }
  \caption{Neighbors of $P$ obtained in Example \ref{Exm.29}
    for a group $I4_122$.} 
\label{Fig.29facets}
\end{figure}

Specially good is the group $I4_122$, whose graphical representation
appears in part (a) of Figure \ref{Fig.29facets}. It has screw
rotations or order-4 in the two versions `dextro' and `levo', and
allows a point to be considered as lying in two different `helices'.

\begin{Exm}
\label{Exm.29}
The following base point and 
metric parameters for a tetragonal group $I4_122$ produce Dirichlet
stereohedra with 29 facets:
    \begin{center}
    Minimal length of horizontal translation = 4

    Minimal length of vertical translation = 1

    Base point $P=(1, \frac{1}{2}, \frac{1}{16})$, in the coordinate system 
    of Figure \ref{Fig.29facets} (a)
    \end{center}
\end{Exm}
Part (b) of Figure \ref{Fig.29facets} shows the 29 orbit
    points producing facets of the stereohedron, with the 
    same conventions of the previous examples. The zero in one of
    the fundamental subdomains indicates that the corresponding neighbor
    is at the same height as the base point $P$.

\end{document}

%% file: table97.tex
{\small{
\begin{table}[htb]
\begin{center}
\renewcommand{\arraystretch}{1.2} 
\begin{tabular}[t]{|c|cc|c|}
\hline
System & Point  &(aspects)& Crystallographic groups \\
       & group  &         &   without reflections   \\
\hline
{\bf Triclinic}      & $1$             & (1)   & $P1$              \\
\cline{2-4}
                &  $\overline{1}$ & (2)    & $P\overline{1}$   \\
\hline
{\bf Monoclinic}     & $2$           & (2)     & $P2$  $P2_1$  $B2$ \\
\cline{2-4}
                & $m$           & (2)     & $Pb$  $Bb$         \\
\cline{2-4}
                & $\frac{2}{m}$ & (4)     & $P\frac{2}{b}$ $P\frac{2_1}{b}$  $B\frac{2}{b}$ \\
\hline
{\bf Orthorhombic}    & $222$         & (4)     & $P22$  $P2_122$  $P22_12_1$  $P2_12_12_1$ \\
                      &               &         & $C222$ $C2_122$  $F222$  $I222$ $I222^{'}$ \\
\cline{2-4}
                & $2mm$         & (4)     & $P2cc$  $P2_1ca$  $P2cn$  $P2ba$  $P2_1bn$ \\
                &               &         & $P2nn$  $C2cc$  $I2cc$ $F2dd$ $A2ba$ \\
\cline{2-4}                                        
                &$\frac{2}{m}\frac{2}{m}\frac{2}{m}$ & (8)
                                        & $P\frac{2}{n}\frac{2}{n}\frac{2}{n}$  
                                          $P\frac{2}{a}\frac{2}{n}\frac{2_1}{n}$  
                                          $P\frac{2}{n}\frac{2_1}{c}\frac{2_1}{c}$\\
                &             &         & $P\frac{2}{a}\frac{2_1}{c}\frac{2}{c}$ 
                                          $P\frac{2_1}{a}\frac{2_1}{b}\frac{2_1}{c}$ 
                                          $P\frac{2_1}{n}\frac{2}{c}\frac{2_1}{a}$ \\ 
                &             &         & $P\frac{2}{n}\frac{2}{b}\frac{2}{a}$ 
                                          $C\frac{2}{a}\frac{2}{c}\frac{2}{c}$  
                                          $I\frac{2}{a}\frac{2}{c}\frac{2}{c}$  
                                          $F\frac{2}{d}\frac{2}{d}\frac{2}{d}$ \\

\hline
{\bf Trigonal} & $3$ & (3)                      & $P3$ $P3_1$ $R3$ \\ 
\cline{2-4}
        &$\overline{3}$ & (6)             & $P\overline{3}$ $R\overline{3}$ \\
\cline{2-4}
        &$32$ & (6)                       &$P32$ $P3_12$ $P312$ $P3_112$ $R32$ \\
\cline{2-4}
        &$3m$ & (6)                       & $P3c$ $P3_1c$ $R3c$ \\
\cline{2-4}
        &$\overline{3}\frac{2}{m}$ & (12) & $P\overline{3}\frac{2}{c}$ $P\overline{3}1\frac{2}{c}$ 
                                          $R\overline{3}\frac{2}{c}$ \\
\hline
{\bf Hexagonal} & $6$ & (6)                     & $P6$ $P6_1$ $P6_2$ $P6_3$\\ 
\cline{2-4}
        & $622$ & (12)                    & $P622$ $P6_122$ $P6_222$ $P6_322$ \\ 
\cline{2-4}
        & $ 6mm$ & (12)                   & $P6cc$ \\
\hline
{\bf Tetragonal}      & $4$ & (4)               & $P4$ $P4_1$ $P4_2$ $I4$ $I4_1$ \\ 
\cline{2-4}
                & $ \overline{4}$ & (4)   & $P\overline{4}$ $I\overline{4}$ \\
\cline{2-4}
                &$\frac{4}{m}$ & (8)      & $P\frac{4}{n}$ $P\frac{4_2}{n}$ $I\frac{4_1}{g}$ \\

\cline{2-4}
                &$422$ & (8)              & $P422$ $P4_122$ $P4_222$ $P42_12$ \\
                &      &                  & $P4_12_12$ $P4_22_12$ $I422$ $I4_122$ \\ 

\cline{2-4}
                &$4mm$ & (8)              & $P4_2gc$ $P4nc$ $P4cc$ $I4_1cd$ \\

\cline{2-4}
  &$\overline{4}2m$ & (8)   
   & $P\overline{4}2c$ $P\overline{4}2_1c$ $P\overline{4}g2$ \\
 &&& $P\overline{4}c2$ $P\overline{4}n2$ $I\overline{4}c2$ $I\overline{4}2d$\\

\cline{2-4}
                & $\frac{4}{m}\frac{2}{m}\frac{2}{m}$ & (16)
                                        & $P\frac{4}{n}\frac{2_1}{c}\frac{2}{c}$ 
                                          $P\frac{4_2}{n}\frac{2}{g}\frac{2}{c}$
                                          $P\frac{4}{n}\frac{2}{n}\frac{2}{c}$ 
                                          $I\frac{4_1}{g}\frac{2}{c}\frac{2}{d}$ \\
\hline
\end{tabular}
\caption{The 97 non-cubic 3-dimensional crystallographic groups 
without reflections}
\renewcommand{\arraystretch}{1} 
\label{Table.Groups}
\end{center}
\end{table}

}}

%% file: table58.tex
\begin{table}  
 \renewcommand{\arraystretch}{1.2}  
\begin{tabular}{|c||c|c||ccc|c||cc|} 
\hline 
Group    & $a$ & $G_0$ & $a_0$ & $i$ & $l$ & Cor. \ref{Coro.planes}
& Final bound &Remarks \\ 
\hline 
 $P\frac{2}{n}\frac{2}{n}\frac{2}{n}$       
                & 8 & $p2$  & 2 & 7 & 1 & 50  & 38 & Prop. \ref{Prop.rectangular.p2} \\
 $P\frac{2}{a}\frac{2}{n}\frac{2_1}{n}$     
                & 8 & $p2$  & 2 & 7 & 1 & 50  & 38 & Prop. \ref{Prop.rectangular.p2} \\
 $P\frac{2}{n}\frac{2_1}{c}\frac{2_1}{c}$   
                & 8 & $p2$  & 2 & 7 & 1 & 50  & 38 & Prop. \ref{Prop.rectangular.p2} \\
 $P\frac{2}{a}\frac{2_1}{c}\frac{2}{c}$     
                & 8 & $p2$  & 2 & 7 & 1 & 50  & 38 & Prop. \ref{Prop.rectangular.p2} \\
 $P\frac{2_1}{a}\frac{2_1}{b}\frac{2_1}{c}$ 
                & 8 & $pg$  & 2 & 7 & 1 & 50  & 38   & Prop. \ref{Prop.pg} \\
 $P\frac{2_1}{n}\frac{2}{c}\frac{2_1}{a}$   
                & 8 & $pg$  & 2 & 7 & 1 & 50  & 38   & Prop. \ref{Prop.pg} \\
 $P\frac{2}{n}\frac{2}{b}\frac{2}{a}$       
                & 8 & $pgg$ & 4 & 7 & 1 & 22  &--- & --- ---\\
 $C\frac{2}{a}\frac{2}{c}\frac{2}{c}$       
                & 8 & $p2$  & 2 & 7 & 1 & {\bf 50}  & ---   & --- ---\\
 $I\frac{2}{a}\frac{2}{c}\frac{2}{c}$       
                & 8 & $pgg$ & 4 & 7 & 2 & {\bf 50}  &  ---  & --- ---  \\
 $F\frac{2}{d}\frac{2}{d}\frac{2}{d}$       
                & 8 & $p2$  & 2 & 7 & 2 &106  & {\bf 70} & Delone \\ 
\hline  

\end{tabular} 
 
 \renewcommand{\arraystretch}{1}  
\caption{Orthorhombic groups without reflexions with more than 4
aspects}
\label{Table.Orthorhombic} 
\end{table}

\begin{table}  
 \renewcommand{\arraystretch}{1.2}  
\begin{tabular}{|c||c|c||ccc|c||cc|}
\hline
Group    & $a$ & $G_0$ & $a_0$ & $i$ & $l$ & Cor. \ref{Coro.planes}
& Final bound &Remarks \\
\hline 
$P\overline{3}$ & 6 & $p3$ & 3 & 4 & 1 & 16 & ---& --- ---\\
$R\overline{3}$ & 6 & $p3$ & 3 & 4 & 1 & 48 & {\bf 42} & Prop. \ref{Prop.R3}\\
\hline
$P32$     & 6 & $p3$ & 3 & 4 & 1 & 16 & ---& --- ---\\
$P3_12$   & 6 & $p1$ & 1 & 4 & 1 & {\bf 48} & ---& --- ---\\
$P312$    & 6 & $p3$ & 3 & 4 & 1 & 16 & ---& --- ---\\
$P3_112$  & 6 & $p1$ & 1 & 4 & 1 & {\bf 48} & ---& --- ---\\
$R32$     & 6 & $p3$ & 3 & 4 & 3 & 48 & {\bf 42} & Prop. \ref{Prop.R3}\\
\hline
$P3c$    & 6 & $p3$ & 3 & 4 & 1 &16 & ---& --- ---\\
$P31c$   & 6 & $p3$ & 3 & 4 & 1 &16 & ---& --- ---\\
$R3c$    & 6 & $p3$ & 3 & 4 & 3 & 48 & {\bf 42} & Prop. \ref{Prop.R3}\\
\hline
$P \overline{3} \frac{2}{c}$  & 12 & $p3$ & 3 & 4 & 1 & 32 & ---&  --- ---\\
$P \overline{3}1 \frac{2}{c}$ & 12 & $p3$ & 3 & 4 & 1 & 32 & ---&  --- ---\\
$R \overline{3} \frac{2}{c}$  & 12 & $p3$ & 3 & 4 & 3 & 96 & {\bf 79} 
                                                     & Cor. \ref{Coro.R32c}\\ 
\hline
\end{tabular} 
 
 \renewcommand{\arraystretch}{1}  
\caption{Trigonal groups without reflexions with more than 4
aspects}
\label{Table.Trigonal} 
\end{table}

\begin{table}  
 \renewcommand{\arraystretch}{1.2}  
\begin{tabular}{|c||c|c||ccc|c||cc|}
\hline
Group    & $a$ & $G_0$ & $a_0$ & $i$ & $l$ & Cor. \ref{Coro.planes}
& Final bound &Remarks \\
\hline 
$P6$   & 6 & $p6$ & 6 & 4 & 1 & 8 & --- & --- ---\\
$P6_1$ & 6 & $p1$ & 1 & 4 & 1 & {\bf 48} & ---& --- ---\\
$P6_2$ & 6 & $p2$ & 2 & 7 & 1 & 36 & --- & --- ---\\
$P6_3$ & 6 & $p3$ & 3 & 4 & 1 & 16 & ---  &--- ---\\
\hline
$P 6 c c$   & 12 & $p6$ & 6 & 4 & 1 & 16 & ---&\\ 
\hline
$P 6 2 2$   & 12 & $p6$ & 6 & 4 & 1 & 16 &---&--- ---  \\ 
$P 6_3 2 2$ & 12 & $p3$ & 3 & 4 & 1 & 32 &---&--- ---   \\ 
$P 6_1 2 2$ & 12 & $p1$ & 1 & 4 & 1 & 96 & {\bf 78} &Prop. \ref{Prop.P6122} \\ 
$P 6_2 2 2$ & 12 & $p2$ & 2 & 7 & 1 & {\bf 78} &---&--- --- \\ 
\hline

\end{tabular} 
 
 \renewcommand{\arraystretch}{1}  
\caption{Hexagonal groups without reflexions with more than 4
aspects}
\label{Table.Hexagonal} 
\end{table}

\begin{table}  
 \renewcommand{\arraystretch}{1.2}  
\begin{tabular}{|c||c|c||ccc|c||cc|}
\hline
Group    & $a$ & $G_0$ & $a_0$ & $i$ & $l$ & Cor. \ref{Coro.planes}
& Final bound &Remarks \\
\hline
$P\frac{4}{n}$   & 8 & $p4$ & 4 & 4 & 1 & 16 & ---&--- ---\\  
$P\frac{4_2}{n}$ & 8 & $p2$ & 2 & 7 & 1 & 50 & 38 & Prop. \ref{Prop.rectangular.p2}\\  
$I\frac{4_1}{g}$ & 8 & $p2$ & 2 & 7 & 1 & 106& {\bf 70}& Delone\\  
\hline
$P422$      & 8  & $p4$  & 4 & 4 & 1 & 16 & ---&--- ---\\ 
$P4_122$    & 8  & $p1$  & 1 & 4 & 1 & 64 & {\bf 50} & Cor. \ref{Coro.P4122}  \\ 
$P4_222$    & 8  & $p2$  & 2 & 7 & 1 & 50 & 38 & Prop. \ref{Prop.rectangular.p2}   \\ 
$P42_12$    & 8  & $p4$  & 4 & 4 & 1 & 16 & ---&--- ---\\ 
$P4_12_12$  & 8  & $p1$  & 1 & 4 & 1 & {\bf 64} & ---&  --- ---  \\ 
$P4_22_12$  & 8  & $p2$  & 2 & 7 & 1 & 50 & 38 & Prop. \ref{Prop.rectangular.p2}\\ 
$I422$      & 8  & $p4$  & 4 & 4 & 2 & 32 & ---&--- ---\\ 
$I4_122$    & 8  & $p2$  & 2 & 7 & 2 &106 & {\bf 70}   & Delone  \\ 
\hline
$P4_2gc$  & 8  & $pgg$ & 4 & 7 & 1 & 22 & ---&--- ---\\  
$P4nc$    & 8  & $p4$  & 4 & 4 & 1 & 16 & ---&--- ---\\  
$P4cc$    & 8  & $p4$  & 4 & 4 & 1 & 16 & ---&--- ---\\  
$I4_1cd$  & 8  & $pgg$ & 4 & 7 & 2 & 50 & {\bf 44} & Prop. \ref{Prop.pgg} \\  
\hline
$P\overline{4}2c$   & 8 & $p2$  & 2 & 7 & 1 & 50 & 38 & Prop. \ref{Prop.rectangular.p2}\\
$P\overline{4}2_1c$ & 8 & $p2$  & 2 & 7 & 1 & 50 & 38 & Prop. \ref{Prop.rectangular.p2}\\
$P\overline{4}g2$   & 8 & $pgg$ & 4 & 7 & 1 & 22 & ---&--- ---\\
$P\overline{4}c2$   & 8 & $p2$  & 2 & 7 & 1 & 50 & 38 & Prop. \ref{Prop.rectangular.p2}\\
$P\overline{4}n2$   & 8 & $p2$  & 2 & 7 & 1 & 50 & 38 & Prop. \ref{Prop.rectangular.p2}\\
$I\overline{4}c2$   & 8 & $pgg$ & 4 & 7 & 2 & 50 & {\bf 40} & Cor. \ref{Coro.pgg.vert} \\
$I\overline{4}2d$   & 8 & $p2$  & 2 & 7 & 2 &106 &{\bf 70} & Delone  \\
\hline 
$P\frac{4}{n}\frac{2_1}{c}\frac{2}{c}$ & 16 & $p4$  & 4 & 4 & 1 & 32 & ---&--- ---\\ 
$P\frac{4}{n}\frac{2}{n}\frac{2}{c}$   & 16 & $p4$  & 4 & 4 & 1 & 32 & ---&--- ---\\ 
$P\frac{4_2}{n}\frac{2}{g}\frac{2}{c}$ & 16 & $pgg$ & 4 & 7 & 1 & 50 & {\bf 40} & Cor. \ref{Coro.pgg.vert} \\
$I\frac{4_1}{g}\frac{2}{c}\frac{2}{d}$ & 16 & $pgg$ & 4 & 7 & 2 &106 &{\bf 80} & Cor. \ref{Coro.pgg.vert} \\
 
\hline 
 
\end{tabular} 
 
 \renewcommand{\arraystretch}{1}  
\caption{Tetragonal groups without reflexions with more than 4
aspects}
\label{Table.Tetragonal} 
\end{table}